\documentclass[a4paper,preprint,11pt,authoryear]{elsarticle}
\usepackage{geometry}
\usepackage{amssymb}
\usepackage{subfigure}
\usepackage{amsmath}
\usepackage{booktabs}
\usepackage{multirow}
\usepackage{array}
\usepackage{lineno}
\usepackage{algorithm}
\usepackage{algpseudocode}
\usepackage{threeparttable}
\usepackage{amsthm}
\usepackage{makecell}
\usepackage{setspace}

\theoremstyle{remark}

\journal{EUROPEAN JOURNAL OF OPERATIONAL RESEARCH}
\allowdisplaybreaks
\begin{document}

\begin{frontmatter}



\title{The Dynamic Team Orienteering Problem in Spatial Crowdsourcing: A Scenario Sampling Approach} 

\affiliation[SCU]{organization={Business School, Sichuan University},
	city={Chengdu},
	postcode={610065},
	country={China}}
	
\author[SCU]{Zhibin Wu\corref{cor1}}\ead{zhibinwu@scu.edu.cn}
\author[SCU]{Songhao Shen}\ead{ssh_nebula@163.com}
\author[SCU]{Yufeng Zhou}
\author[SCU]{Qin Lei}
\cortext[cor1]{Corresponding author}


%

\begin{abstract}
	In services such as retail audits and urban infrastructure monitoring, a platform dispatches rewarded, location-based micro-tasks to mobile workers traveling along personal origin-destination (OD) trips under hard time budgets. 
	As requests with time constraints arrive online over a finite horizon, the platform must decide which requests to accept and how to route workers to maximize collected profit. 
	We model this setting as the Dynamic Team Orienteering Problem in Spatial Crowdsourcing (DTOP-SC).
	To solve this problem, we propose a scenario-sampling rolling-horizon framework that mitigates myopic bias by augmenting each planning epoch with sampled virtual tasks. 
	At each epoch, the augmented task set defines a deterministic static subproblem solved via an adaptive large neighborhood search (ALNS).
	We also formulate a mixed-integer programming model to provide offline reference solutions.
	Computational experiments are conducted on synthetic DTOP-SC instances generated from real-world road-map coordinates and on a dynamic team orienteering (DTOP) benchmark.
	On the map-based instances, the proposed policy exhibits stable gaps with respect to time-limited MIP solutions across the tested scales, while maintaining smooth computational scalability as the problem size increases. 
	On the DTOP benchmark, the policy achieves an average decision time of 0.14s per instance, with 192-198s reported for multiple plan approach as an indicative reference, while maintaining competitive profit.
\end{abstract}



\end{frontmatter}


\section{Introduction}

\begin{figure*}[!tbp]
	\centering
	\includegraphics[width=0.65\linewidth]{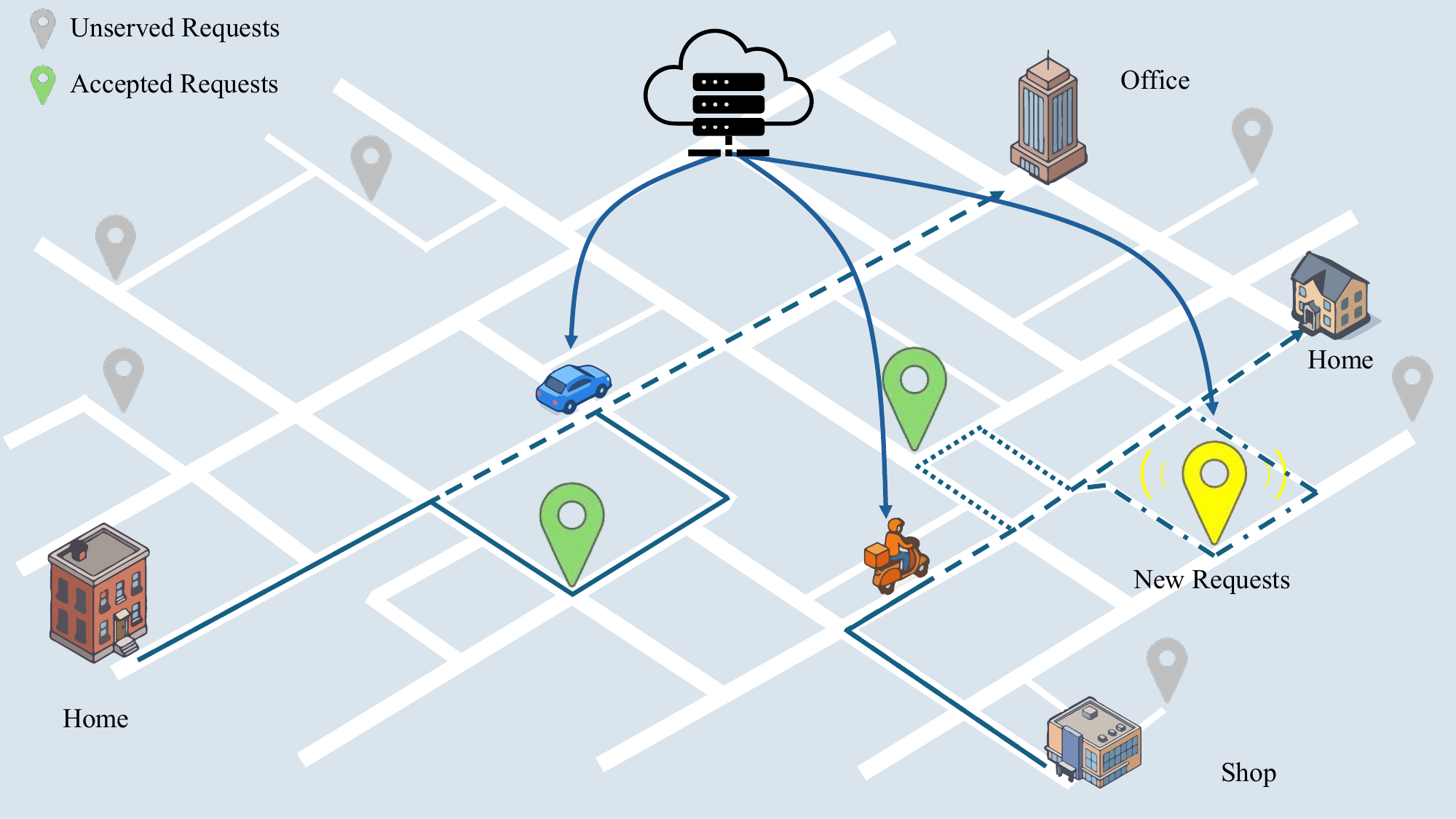}
	\caption{Schematic illustration of the Dynamic Team Orienteering Problem in Spatial Crowdsourcing.}
	\label{fig:problem}
\end{figure*}

Over the past decade, a class of on-demand services has seen growing adoption, in which a platform uses mobile workers to collect timely, fine-grained information about the physical world. 
Examples include retail audit services in brick-and-mortar stores, urban infrastructure monitoring, and geo-tagged observations of traffic, environment, or local events. 
Spatial crowdsourcing and mobile crowdsensing platforms support such services by leveraging mobile device carriers as roaming sensors and assigning location-based micro-tasks to them \citep{tong2020,suhag2023}. 
In these settings, a platform dispatches geographically distributed information tasks with rewards to a pool of roaming workers. 
Workers briefly visit task locations to take photos, record measurements, or verify local conditions for remote requesters.
Workers are often part-time participants, such as commuters or everyday shoppers, who are only willing to accept tasks that can be bundled into their personal trips and completed within individual time budgets. 
As illustrated in Figure~\ref{fig:problem}, the platform assigns online tasks to these mobile workers traveling between specific origins and destinations (e.g., Home to Office). 
Because not all tasks can be served, the platform must jointly decide task selection, worker assignment, and service sequencing. 
The objective is to maximize the total collected information value subject to feasibility constraints.

Two structural features of these settings are particularly salient. 
First, the environment is dynamic: new tasks, service requirements, and travel conditions are progressively revealed over time, so that assignment and routing decisions must be revised repeatedly under partial information. 
Second, the workforce is heterogeneous in both trajectories and temporal constraints: workers have different origins/destinations and worker-specific availability windows and time budgets. 
This heterogeneity yields worker-dependent feasible regions and substantially different opportunity costs for accepting new tasks.
When requests arrive online, ignoring heterogeneous OD trajectories and remaining working times can lead to assigning distant tasks to poorly positioned workers. 
This may cause end-of-shift infeasibilities and long empty repositioning, while leaving later tasks along other workers' low-detour corridors unserved. 
Hence, dynamic assignment and routing decisions need to account explicitly for how each near-term commitment interacts with a worker's remaining time budget and residual travel path.

Motivated by these applications, we study a dynamic profitable routing problem. 
A central decision maker controls a heterogeneous team of workers to serve rewarded location-based information tasks that arrive online over a finite planning horizon. 
Each task specifies a service location, a service duration, a time window and release time, and yields a profit if served. 
Each worker has an individual OD pair and a personal working-time window, within which they must complete all assigned tasks and still reach their destination. 
At successive decision epochs, the platform decides which tasks to accept, assigns accepted tasks to workers, and determines the within-shift service sequence of each worker based on the information revealed so far. 
The objective is to maximize the total collected profit while satisfying all temporal and trajectory constraints. 
We refer to this problem as the Dynamic Team Orienteering Problem in Spatial Crowdsourcing (DTOP-SC) to reflect its spatial crowdsourcing setting.
From a structural and taxonomic perspective in orienteering, the DTOP-SC falls under the Dynamic Heterogeneous-Trajectory Team Orienteering Problem with Time Windows (D-HT-TOPTW), characterized by worker-specific OD trips with hard time budgets and tasks with release times and time windows.
Hereafter, we use DTOP-SC to denote this problem.
While we use a common travel-time matrix in our notation and experiments, the formulation and solution framework can accommodate worker-specific travel times by allowing worker-dependent travel-time matrices.

From a modeling perspective, DTOP-SC extends TOPTW to an online setting with worker-specific OD pairs and working-time windows, capturing heterogeneous feasible corridors and opportunity costs in spatial information services. 
When all tasks are known in advance, it reduces to a static counterpart, Heterogeneous-Trajectory Team Orienteering Problem with Time Windows (HT-TOPTW), for which we provide a mixed-integer programming (MIP) formulation as a time-limited offline baseline. 
DTOP-SC also subsumes the multi-start multi-end TOPTW of \citet{shen2025} and the DTOP of \citet{kirac2025} as special cases.
To mitigate the myopia of purely reactive rolling-horizon dispatching, we develop Scen-RH-ALNS, a scenario-sampling rolling-horizon ALNS framework (abbreviated as S-ALNS in tables for brevity).
At each decision epoch, we augment the current snapshot with sampled virtual tasks and solve the resulting deterministic instances in parallel using ALNS. 
We then aggregate robust near-term assignments across scenarios to form a joint dispatch decision.

The main contributions of this paper are threefold.

(1) We formulate DTOP-SC, a dynamic heterogeneous trajectory team orienteering problem with time-windowed tasks in spatial crowdsourcing. 
We also provide a MIP formulation for its static counterpart as a time-limited offline baseline.

(2) A parallelizable scenario-sampling rolling-horizon ALNS designed for real-time task assignment and routing in DTOP-SC is proposed. 
By augmenting static snapshots with sampled virtual tasks, this framework provides an explicit lookahead mechanism to mitigate myopic bias in dynamic dispatching.

(3) Extensive computational experiments are conducted on synthetic DTOP-SC instances based on real-world road-map coordinates and on standard DTOP benchmarks. 
The results show that Scen-RH-ALNS achieves routing quality comparable to state-of-the-art multiple plan approaches in medium-to-high dynamism settings. 
It also delivers sub-second decision times on DTOP benchmarks and requires only a few seconds on larger map-based instances.
These runtimes correspond to two to three orders of magnitude shorter elapsed times than the reported runtimes of existing multi-plan approaches.

The remainder of this paper is organized as follows. 
Section~\ref{sec:related} reviews related work on the team orienteering problems and  dynamic routing methods. 
Section~\ref{sec:problem} defines the dynamic DTOP-SC problems and presents the MIP formulation for its static counterpart. 
Section~\ref{sec:method} describes the proposed Scen-RH-ALNS framework in detail. 
Section~\ref{sec:experiments} presents the experimental setup and computational results on benchmark and map-based instances. 
Finally, Section~\ref{sec:conclusion} concludes and discusses directions for future research.

\section{Related work}
\label{sec:related}

We briefly review two strands of literature related to our work: Team orienteering problems and dynamic routing methods.

\subsection{Team orienteering problems}

Our modeling framework builds on the literature on orienteering and selective-routing problems. 
In the Orienteering Problem (OP), a subset of locations is visited to maximize the total collected profit subject to a travel-time or distance budget \citep{tsiligirides1984,golden1987,vansteenwegen2011}. 
This basic model has been extended in many directions, including multiple routes (the Team Orienteering Problem, TOP), time windows (TOPTW), and numerous application-driven variants \citep{chao1996,vansteenwegen2009,gunawan2016}. 
These models retain the profit-maximizing, selective-routing core and have been applied in contexts such as tourist trip design, maintenance planning, and more recently crowdsensing and spatial crowdsourcing \citep[e.g.,][]{gunawan2017,gedik2017,hammami2020,chaigneau2025,perez2025}.

Several extensions relax the single-depot and closed-tour assumptions, for example, multi-depot open-route formulations in industrial contexts \citep[e.g.,][]{chaigneau2025}, and the MSMETOPTW of \citet{shen2025} in spatial crowdsourcing, where each worker has an individual origin, destination, and working-time budget. 
Dynamic and stochastic generalizations have also been proposed. 
Dynamic variants reveal customers online over a fixed planning horizon \citep{kirac2025}, while stochastic/probabilistic variants model random request arrivals or release dates and are often formulated as MDPs \citep[e.g.,][]{angelelli2021,li2024}; learning-based policies have also been explored under dynamic/uncertain travel conditions \citep{ammouriova2024}. 
These contributions demonstrate the versatility of orienteering models in capturing online arrivals and uncertainty, and they provide important methodological building blocks for our work.
Much of the existing literature focuses on depot-based operations and largely homogeneous fleets, or abstracts worker availability without explicitly modeling individual OD trips.
This motivates studying OD-based crowdsourcing settings in which feasibility corridors and opportunity costs are worker-specific.

From this perspective, the problem studied in this paper (DTOP-SC) corresponds to a dynamic heterogeneous-trajectory TOPTW variant (D-HT-TOPTW), where workers follow individual OD trips under worker-specific working-time windows. 
It subsumes the MSMETOPTW of \citet{shen2025} and the depot-based DTOP of \citet{kirac2025} as special cases.
We also introduce a static counterpart (HT-TOPTW) and provide a MIP formulation to serve as a time-limited offline baseline.

\subsection{Dynamic routing methods}
Beyond the orienteering family, our work is related to dynamic routing and dynamic vehicle routing problems (DVRPs). 
Surveys such as \citet{pillac2013}, \citet{psaraftis2016}, and \citet{ritzinger2016} review the dynamic vehicle routing literature. 
In particular, \citet{pillac2013} and \citet{psaraftis2016} classify dynamic routing problems based on how information unfolds. 
In dynamic-deterministic settings, future demand is unknown and no explicit probabilistic model is assumed. 
In dynamic-stochastic settings, forecasts or probability distributions are exploited.
The DTOP-SC studied here can be viewed as a dynamic-deterministic, profit-maximizing routing problem with time windows and a heterogeneous OD-based workforce.

On the algorithmic side, dynamic routing is often addressed via a rolling-horizon framework \citep{psaraftis2016}. 
At each decision epoch, the dispatcher freezes the current state and solves a static approximation. 
The solution is implemented only partially (typically the first route segment), and the process repeats when new information becomes available. 
Static subproblems are often solved by metaheuristics such as (A)LNS. 
ALNS has been successfully applied to rich VRP variants \citep{ropke2006} and to TOP \citep{hammami2020}, while LNS variants are also competitive for very-large-scale TOP \citep{chaigneau2025}. 
The rolling-horizon framework has been deployed in dynamic transit, urban parcel logistics, and multi-period delivery systems with endogenous demand \citep[e.g.,][]{ma2023,faugere2022,keskin2023}. 
As discussed in \citet{ma2023}, a rolling-horizon approach that only uses current information can be short-sighted because it overlooks the impact of future requests on routing decisions. In our profit-maximizing setting, this limitation also relates to the opportunity cost of committing limited time budgets to currently known tasks.

Anticipatory policies address this limitation by exploiting information about future demand, for example through scenario-based planning. 
In the seminal framework of \citet{bent2004}, sampled demand scenarios are solved as static VRPs and near-term actions are extracted from these solutions, for example via consensus function. 
Variants of this multiple-scenario paradigm have been applied to dynamic time-window assignment in service routing, explicitly trading off immediate routing cost against preserving capacity for future requests \citep{paradiso2025}. 
In dynamic orienteering, \citet{kirac2025} proposed a multiple plan approach for the DTOP, which maintains a diverse pool of alternative routing plans and selects actions online.
This strategy provides robustness to uncertainty without committing to an explicit stochastic demand model.
In multi-worker settings with richer constraints, maintaining and improving a sufficiently diverse plan pool may become computationally demanding, motivating complementary anticipatory schemes with lightweight scenario generation.

Our work follows the rolling-horizon paradigm and employs ALNS to solve static snapshots, but incorporates scenario sampling at each decision epoch to construct virtual tasks representing possible future requests. 
By extracting and aggregating robust near-term actions across scenarios, we obtain a multi-worker joint dispatch policy with explicit lookahead that remains computationally tractable. 
Accordingly, our approach bridges purely reactive rolling-horizon re-optimization and plan-diversification strategies by introducing an explicit yet lightweight lookahead mechanism via scenario sampling, while retaining real-time tractability in heterogeneous multi-worker settings.

\section{Problem formulation}
\label{sec:problem}

\subsection{Operational setting}
\label{subsec:operational_setting}

We consider an on-demand service in which a central platform dispatches mobile workers to visit spatially distributed tasks and collect in-situ information on behalf of remote requesters. 
Let $\mathcal{W}$ denote the finite set of workers and $V_{\mathcal{T}}$ the finite set of tasks that may potentially appear during the planning horizon $[0,H]$. 
Each worker $w \in \mathcal{W}$ has an origin $s_w$ and a destination $d_w$, interpreted as their home or workplace, and a personal working-time window $[T_{start}^w,T_{end}^w]$ within which all assigned tasks must be completed and the worker must reach $d_w$. 
Workers are assumed to travel in a transportation environment characterized by a given travel-time matrix.
Let $V = V_{\mathcal{T}} \cup \{s_w,d_w\}_{w \in \mathcal{W}}$ denote the set of all relevant nodes (tasks, worker origins, and worker destinations) in the transportation environment. 
For notational simplicity, we assume a common deterministic travel-time matrix $(t_{ij})_{i,j \in V}$ shared by all workers, where $t_{ij} \ge 0$ denotes the travel time from node $i$ to node $j$.
All formulations and algorithms in this paper extend directly to worker-specific travel times $t_{ij}^w$. 
In this way, we can model heterogeneous travel speeds or transportation modes across workers. 
We maintain the assumptions that travel times are deterministic and do not depend on departure time.

Each task $i \in V_{\mathcal{T}}$ corresponds to visiting a physical location once to perform a short information-collection operation, such as taking photos, recording measurements, or verifying local conditions. 
Task $i$ yields a profit $p_i \ge 0$ to the platform if it is served. 
The service at $i$ requires a deterministic duration $\tau_i \ge 0$ and must start within a time window $[b_i,e_i] \subseteq [0,H]$. 
In addition, each task has a release time $r_i \in [0,H]$ at which it becomes known to the dispatcher. 
We assume that $r_i \le b_i \le e_i$ for all tasks, so that the platform never learns about a task after its time window has closed. 
If task $i$ is not served by any worker within its time window, it expires without generating profit.

The platform has full knowledge of the travel-time matrix, worker origins and destinations, working-time windows, and the distributions or empirical ranges of task attributes. 
However, task arrivals are dynamic: only tasks with $r_i \le t$ are known at time $t$, and the platform must make routing and assignment decisions without observing future tasks. 
Travel times and service durations are assumed deterministic and time-independent; congestion and stochastic travel times are beyond the scope of this work.

\subsection{Dynamic team orienteering problem in spatial crowdsourcing}
\label{subsec:dynamic_DH-TOPTW}

We now formalize the resulting dynamic routing problem, which we call the Dynamic Team Orienteering Problem in Spatial Crowdsourcing. 
Time evolves over the horizon $[0,H]$, and tasks appear dynamically at their release times. 
The platform observes task releases over time as well as the evolving states of all workers. 
At each decision epoch, it decides which tasks to accept and assign. 
It then determines the service sequence, aiming to maximize the total collected profit under temporal and trajectory constraints.

At any time $t \in [0,H]$, the set of tasks that have appeared and available for assignment is
\[
A(t) = \{\, i \in V_{\mathcal{T}} : r_i \le t \le e_i,\ i \text{ available for assignment}\,\}.
\]
Workers may be travelling between nodes, serving a task, or idle and available for new assignments. 
For each worker $w \in \mathcal{W}$, we maintain a planned route
\[
\pi_w^t = (s_w, i_1, i_2, \dots, i_{g-1}, i_g, \dots, d_w),
\]
consisting of its origin $s_w$, a sequence of intermediate task nodes, and its destination $d_w$. 
The executed prefix $R_w^t$ of $\pi_w^t$ contains all nodes that $w$ has already visited by time $t$, and the remaining suffix $\pi_w^t \setminus R_w^t$ contains future planned visits, which may be revised as new information arrives. 
Let $\operatorname{pred}(w,t)$ denote the last node visited by $w$ before or at time $t$ and $\operatorname{succ}(w,t)$ the next planned node on $\pi_w^t$; then
\[
R_w^t = (s_w, i_1, \dots, \operatorname{pred}(w,t)), \qquad
\pi_w^t \setminus R_w^t = (\operatorname{succ}(w,t), \dots, d_w).
\]

Service start times along a planned route are computed recursively. 
For a route
$\pi_w^t = (i_0,i_1,\dots,i_G)$ with $i_0 = s_w$ and $i_G = d_w$, we set
\[
a_0^w = T_{start}^w, \qquad
a_g^w = a_{g-1}^w + \tau_{i_{g-1}} + t_{i_{g-1},i_g} + \omega_g^w,
\quad g = 1,\dots,G,
\]
where $\omega_g^w \ge 0$ is the waiting time at node $i_g$. 
Feasibility requires that, for every task node $i_g \in V_{\mathcal{T}}$ on the route of worker $w$, the start of service respects the time window and release time,
\[
a_g^w \ge \max\{b_{i_g}, r_{i_g}\}, \qquad a_g^w \le e_{i_g},
\]
and that arrival at the destination satisfies $a_G^w \le T_{end}^w$.

At each time $t$, we define the set of idle workers as
\[
\mathcal{W}_{\mathrm{idle}}(t)
= \{\, w \in \mathcal{W} : T_{start}^w \le t,\ \text{$w$ has not reached $d_w$ and
	is not currently travelling or serving}\,\}.
\]
Workers with $T_{start}^w > t$ are not yet available and are therefore excluded from $\mathcal{W}_{\mathrm{idle}}(t)$.
The dynamic state of the system at time $t$ can then be represented by
\[
X(t) = \bigl(A(t), \{R_w^t,\pi_w^t\}_{w \in \mathcal{W}}\bigr),
\]
which encodes the available tasks, the executed prefixes, and the current planned routes of all workers.

Decisions are taken at a sequence of decision epochs
$0 = \mathbb{T}_0 < \mathbb{T}_1 < \mathbb{T}_2 < \dots \le H$, triggered by operational events such as task arrivals and workers completing service and becoming idle. 
At a generic epoch $t = \mathbb{T}_k$, the platform observes the current state $X(t)$ and chooses a dispatch decision that (i) selects a subset of tasks from $A(t)$ to accept, and (ii) assigns some of these tasks to idle workers as their next destinations, implicitly updating the future planned routes $\{\pi_w^t\}_{w \in \mathcal{W}_{\mathrm{idle}}(t)}$. 
Between decision epochs, workers follow their most recently assigned routes without further re-optimization.

A feasible policy is a rule that specifies, at every decision epoch and for every admissible state $X(t)$, which tasks to accept and which worker-task assignments to implement, subject to the temporal and trajectory constraints introduced above. 
The objective of DTOP-SC is to find a feasible policy that maximizes the total collected profit over the horizon,
\[
\sum_{i \in V_{\mathcal{T}}} p_i \cdot \mathbf{1}\{i \text{ is served within} [b_i,e_i]\},
\]
where $\mathbf{1}\{\cdot\}$ denotes the indicator function. 
In the next subsection, we show that, at any decision epoch, freezing the current state and ignoring future arrivals leads to a static snapshot problem that can be modeled as a Heterogeneous Trajectory Team Orienteering Problem with Time Windows (HT-TOPTW).

\subsection{Static HT-TOPTW snapshot and MIP formulation}
\label{subsec:ht_toptw_mip}

We then present a MIP formulation for the static HT-TOPTW, which arises as a snapshot of DTOP-SC when all tasks are assumed to be known in advance and no new tasks arrive in the future.

Let $V_S = \{s_w\}_{w \in \mathcal{W}}$ and $V_D = \{d_w\}_{w \in \mathcal{W}}$ denote the sets of worker-specific origins and destinations, respectively. 
For notational convenience we define the global node set $V = V_{\mathcal{T}} \cup V_S \cup V_D$.
However, each worker $w$ can only traverse its own OD pair.
Hence we define a worker-specific node set $V^w = V_{\mathcal{T}} \cup \{s_w,d_w\}$, and we only define routing and timing variables on $V^w$.
In particular, $x_{ij}^w$ is defined only for $i,j \in V^w$, and $a_i^w$ is defined only for $i \in V^w$.
We also set $\tau_{s_w}=\tau_{d_w}=0$ for all $w$.

Each task $i \in V_{\mathcal{T}}$ has a profit $p_i$, service duration $\tau_i$, time window $[b_i, e_i]$, and a release time $r_i$ at which it is revealed to the dispatcher. 
In a purely static HT-TOPTW instance we may set $r_i = 0$ for all $i$, but keeping $r_i$ explicit facilitates the connection between the static model to the dynamic setting. 
Each worker $w \in \mathcal{W}$ has an origin $s_w$, a destination $d_w$, a release time $T_{start}^w$ at $s_w$, and a hard deadline $T_{end}^w$ at $d_w$. 
Using binary variables $x_{ij}^w$ and $y_i^w$, and continuous variables $a_i^w$, the static HT-TOPTW can be formulated as:

\paragraph{Objective function}
\begin{equation}
	\label{eq:obj}
	\max \sum_{i \in V_{\mathcal{T}}} p_i \cdot \Big(\sum_{w \in \mathcal{W}} y_i^w\Big)
\end{equation}

\paragraph{Constraints}
\begin{align}
	\label{eq:start_flow}
	&\sum_{j \in V_{\mathcal{T}} \cup \{d_w\}} x_{s_w, j}^w = 1,
	&& \forall w \in \mathcal{W} \\
	\label{eq:end_flow}
	&\sum_{i \in V_{\mathcal{T}} \cup \{s_w\}} x_{i, d_w}^w = 1,
	&& \forall w \in \mathcal{W} \\
	\label{eq:no_into_start}
	&\sum_{i \in V_{\mathcal{T}} \cup \{d_w\}} x_{i,s_w}^w = 0,
	&&\forall w \in \mathcal{W} \\	
	\label{eq:no_out_of_dest}
	&\sum_{j \in V_{\mathcal{T}} \cup \{s_w\}} x_{d_w,j}^w = 0,
	&&\forall w \in \mathcal{W}	\\	
	\label{eq:flow_conservation}
	&\sum_{i \in V^w \setminus \{k\}} x_{ik}^w
	= \sum_{j \in V^w \setminus \{k\}} x_{kj}^w
	= y_k^w,
	&& \forall k \in V_{\mathcal{T}},\ \forall w \in \mathcal{W} \\
	\label{eq:task_once}
	&\sum_{w \in \mathcal{W}} y_i^w \le 1,
	&& \forall i \in V_{\mathcal{T}} \\
	\label{eq:time_propagation}
	&a_i^w + \tau_i + t_{ij} - a_j^w
	\le \mathbb{M}(1 - x_{ij}^w),
	&& \forall i \in V^w \setminus \{d_w\},\ \forall j \in V^w \setminus \{s_w\},\ i \ne j,\ \forall w \in \mathcal{W} \\
	\label{eq:time_low_bound}
	&a_i^w \ge b_i - \mathbb{M}(1 - y_i^w),
	&& \forall i \in V_{\mathcal{T}},\ \forall w \in \mathcal{W} \\
	\label{eq:time_release}
	&a_i^w \ge r_i - \mathbb{M}(1 - y_i^w),
	&& \forall i \in V_{\mathcal{T}},\ \forall w \in \mathcal{W} \\
	\label{eq:time_up_bound}
	&a_i^w \le e_i + \mathbb{M}(1 - y_i^w),
	&& \forall i \in V_{\mathcal{T}},\ \forall w \in \mathcal{W} \\
	\label{eq:start_time}
	&a_{s_w}^w = T_{start}^w,
	&& \forall w \in \mathcal{W} \\
	\label{eq:deadline}
	&a_{d_w}^w \le T_{end}^w,
	&& \forall w \in \mathcal{W} \\
	\label{eq:binary_vars_x}
	&x_{ij}^w \in \{0, 1\},
	&& \forall i,j \in V^w,\ i\ne j,\ \forall w \in \mathcal{W} \\
	\label{eq:binary_vars_y}
	&y_i^w \in \{0, 1\},
	&& \forall i\in V_{\mathcal{T}},\ \forall w\in\mathcal{W} \\
	\label{eq:continuous_vars}
	&a_i^w \ge 0,
	&& \forall i \in V^w,\ \forall w \in \mathcal{W}
\end{align}
A valid choice is $\mathbb{M} \ge \max_{w\in\mathcal{W}}(T_{\mathrm{end}}^w-T_{\mathrm{start}}^w) + \max_{i\in V_{\mathcal{T}}}\tau_i + \max_{i,j\in V} t_{ij}$, which safely dominates any feasible time difference on an unused arc.

The objective \eqref{eq:obj} maximizes the total profit collected from visited tasks, as in a classical orienteering formulation. 
Constraint sets \eqref{eq:start_flow}-\eqref{eq:flow_conservation} route each worker from its own start node $s_w$ to its destination $d_w$ and enforce flow conservation at task nodes, linking arc variables $x_{ij}^w$ to visit variables $y_i^w$. 
Constraints \eqref{eq:no_into_start}-\eqref{eq:no_out_of_dest} further enforce a path structure by forbidding arcs into $s_w$ and out of $d_w$.
Constraint set \eqref{eq:task_once} ensures that each task is served by at most one worker and thereby couples the worker-specific routes on the shared task set. 
Constraint sets \eqref{eq:time_propagation}-\eqref{eq:time_up_bound} impose temporal feasibility. 
Specifically, \eqref{eq:time_propagation} propagates service start times along used arcs, and \eqref{eq:time_low_bound}-\eqref{eq:time_up_bound} enforce task time-window bounds. 
The additional release-time constraints \eqref{eq:time_release} ensure that service at task $i$ cannot start before $r_i$.
Constraint set \eqref{eq:start_time} fixes departure times at origins, and \eqref{eq:deadline} enforces worker-specific deadlines at destinations.
Constraint sets \eqref{eq:binary_vars_x}-\eqref{eq:continuous_vars} specify the domains of the decision variables.

From the perspective of each worker $w$, constraint sets \eqref{eq:start_flow}-\eqref{eq:deadline} define a feasible orienteering path between $(s_w, d_w)$ under time-window, release-time, and deadline constraints. 
Structurally, this formulation differs from classical depot-based TOPTW in two ways. 
First, each worker $w$ faces an individual orienteering subproblem between $(s_w, d_w)$ with its own time budget and release time, so the feasible task set and temporal flexibility are worker-specific and heterogeneous. 
Second, the coupling across workers is entirely through the shared task set via \eqref{eq:task_once}, so interactions arise only where these heterogeneous feasible regions overlap. 
In the dynamic setting studied later, task release times $r_i$ correspond to appearance times of requests, and the above static formulation serves as an offline reference that clarifies the feasible region of DTOP-SC when all information is known in advance.

The DTOP-SC framework contains several existing orienteering models as special cases. 
In particular, the DTOP of \citet{kirac2025} is recovered under a homogeneous, single-depot specification.

Consider a DTOP-SC instance in which all workers share a common origin and destination, denoted by $v_0$, and have identical working horizons:
\[
s_w = v_0,\quad d_w = v_0,\quad
T_{start}^w = T_{\mathrm{start}},\quad
T_{end}^w = T_{\mathrm{end}},
\qquad \forall w \in \mathcal{W}.
\]
All workers therefore operate from the same depot with the same time budget $T_{\mathrm{end}} - T_{\mathrm{start}}$, and share the same travel-time matrix $(t_{ij})_{i,j \in V}$.

Each customer (task) $i \in V_{\mathcal{T}}$ has a profit $p_i$, service duration $\tau_i$, and a time window $[b_i,e_i] \subseteq [T_{\mathrm{start}}, T_{\mathrm{end}}]$ as in a standard TOPTW setting. 
The release time $r_i$ represents the time at which customer $i$ becomes known to the dispatcher: for customers that are known in advance, $r_i = T_{\mathrm{start}}$, whereas for dynamically revealed customers $r_i$ equals their appearance time.
Under this homogeneous single-depot specialization, the DTOP-SC reduces to a dynamic team orienteering problem with a common depot and time budget, which coincides with the DTOP formulation of \citet{kirac2025}.

Throughout the paper, we use a single travel-time matrix $(t_{ij})_{i,j \in V}$ for notational convenience. 
If worker-specific travel times $t_{ij}^w$ are available, all structural properties of HT-TOPTW and the Scen-RH-ALNS framework remain valid. 
Specifically, $t_{ij}$ is replaced by $t_{ij}^w$ in the time-propagation constraints \eqref{eq:time_propagation} and in the recursive computation of service start times along each route. 
The same replacement is applied in the feasibility checks used during snapshot construction. 
The dynamic control logic and the scenario-based lookahead mechanism are unchanged; only the underlying travel-time parameters become worker-dependent.

\section{Problem-solving methodology}
\label{sec:method}

Building on the DTOP-SC model introduced in Section~\ref{sec:problem}, we now propose an online control policy that can be implemented in an operational setting. 
The policy combines a rolling-horizon mechanism, static HT-TOPTW optimization on decision snapshots, and a scenario-sampling lookahead scheme.

The planning horizon is partitioned into a sequence of decision epochs triggered by operational events such as task arrivals and workers becoming idle. 
At each decision epoch, the dispatcher observes the current dynamic state. 
This state includes the available tasks and the current worker routes. 
The dispatcher then formulates a static HT-TOPTW subproblem to capture feasible near-term routing options for idle workers. 
Finally, it approximately solves the subproblem using an ALNS and obtains a set of candidate routes.
A subset of the first moves on these routes is committed as firm worker-task assignments to be executed until the next decision epoch, while workers follow the most recently issued assignments without further re-optimization between epochs. 
To mitigate the myopic nature of such snapshot-based optimization, the ALNS is augmented with sampled demand scenarios and virtual tasks, and a multi-scenario aggregation rule is used to derive a joint dispatch decision for multiple workers.

\begin{figure}[t]
	\centering
	\includegraphics[width=\linewidth]{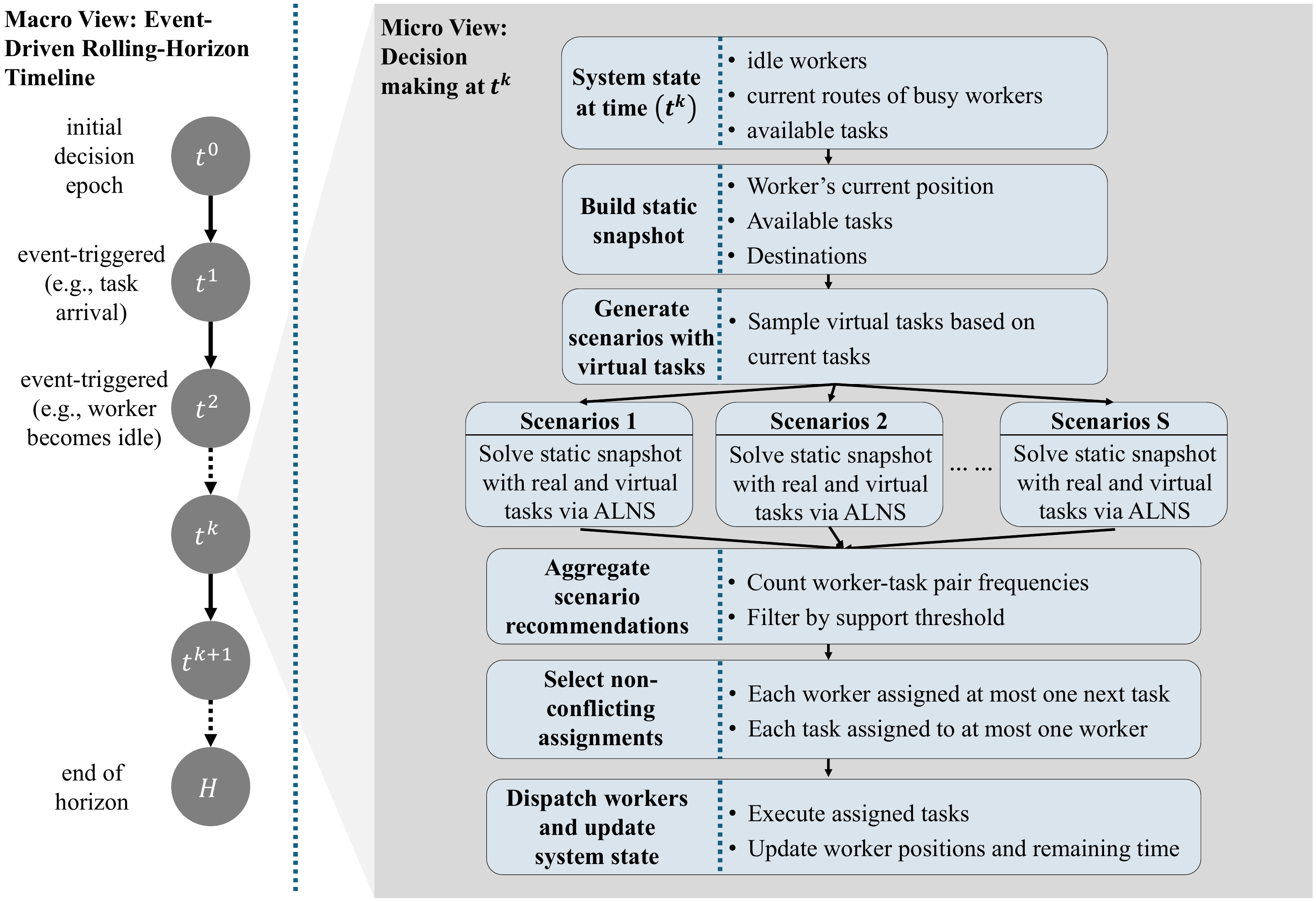}
	\caption{Event-driven rolling-horizon decision process in Scen-RH-ALNS.}
	\label{fig:scen_rh_alns_overview}
\end{figure}

Figure~\ref{fig:scen_rh_alns_overview} schematically illustrates the resulting
event-driven rolling-horizon decision process at both the macro (timeline) and
micro (single decision epoch) levels.

The remainder of this section details the components of the control policy. 
Section~\ref{subsec:rh_alns} presents a myopic rolling-horizon ALNS, which repeatedly constructs and optimizes static HT-TOPTW subproblems based on the currently observed tasks and worker states. 
Section~\ref{subsec:scenario_extension} extends the myopic rolling-horizon ALNS with a scenario-sampling lookahead mechanism and a multi-worker dispatch rule that aggregates candidate actions across scenarios. 
Finally, Section~\ref{subsec:overall} summarizes the overall algorithmic framework and discusses computational aspects.

\subsection{A rolling-horizon ALNS for DTOP-SC}
\label{subsec:rh_alns}

Before introducing the scenario-sampling extension, we first describe a myopic rolling-horizon framework for DTOP-SC. 
This framework provides the static optimization component that will be reused inside the scenario-sampling framework in Section~\ref{subsec:scenario_extension}. 
At a high level, the framework freezes the current dynamic state at each decision epoch. 
It formulates a static HT-TOPTW snapshot over the currently available tasks and idle workers. 
The snapshot is approximately solved using ALNS. The framework then commits a limited set of near-term worker-task assignments.
Only the unexecuted suffixes of the routes of idle workers may be modified; the executed prefixes and the routes of non-idle workers are treated as fixed commitments.

\subsubsection{Initial snapshot construction and static optimization}
\label{subsubsec:initial_snapshot}

At time $t = 0$, a subset $K_0 \subseteq V_{\mathcal{T}}$ of tasks is assumed to be known, typically those with release times $r_i = 0$ or those announced in advance. 
At the same time, the set of workers that are available at time~0, each located at its origin $s_w$, is observed. 
An initial routing plan is constructed by solving a static HT-TOPTW instance defined on the tasks $K_0$ and the workers available at time~0.

Formally, an initial static node set is defined as
\begin{equation}
V(0) = K_0 \cup V_S(0) \cup V_D(0),
\end{equation}
where $V_S(0) = \{s_w\}_{w \in \mathcal{W}_{\mathrm{idle}}(0)}$ and
$V_D(0) = \{d_w\}_{w \in \mathcal{W}_{\mathrm{idle}}(0)}$.
For each worker $w \in \mathcal{W}_{\mathrm{idle}}(0)$, the start and end times used in the snapshot are set to
\[
T_{start}^w(0) = 0, \qquad
T_{end}^w(0) = T_{end}^w,
\]
and for each task $i \in K_0$ the profit $p_i$, service duration $\tau_i$, time window $[b_i,e_i]$, and release time $r_i$ are taken from the original problem data. 
Travel times $t_{ij}$ are inherited from the original instance. 
This yields a static HT-TOPTW instance of the form described by the MIP formulation in Section~\ref{sec:problem}, restricted to tasks $K_0$ and the workers available at time~0. 
Workers whose start times satisfy $T_{start}^w > 0$ do not appear in this initial snapshot and will only enter later snapshots once they become available.

This static HT-TOPTW instance is approximately solved by ALNS. 
The ALNS maintains a set of routes $\{\pi_w^0\}_{w \in \mathcal{W}_{\mathrm{idle}}(0)}$, one per available worker, each starting at $s_w$, ending at $d_w$, and visiting a subset of tasks in $K_0$. 
It starts from a greedy constructive heuristic that builds initial routes by iteratively inserting tasks into worker routes. 
In each iteration, a fraction of the currently routed tasks is removed by one of three destroy operators (Shaw removal, random removal, or worst-cost removal), and the removed tasks are reinserted by one of three repair operators (regret-3 insertion, regret-2 insertion, or greedy insertion). 
The choice of destroy and repair operators follows a roulette-wheel mechanism based on adaptive weights that are periodically updated from their accumulated scores over segments of iterations. 
After each destroy-repair cycle, the solution is further improved by local search: an intra-route 2-opt move is applied systematically, and inter-route relocate and inter-route swap moves are invoked at regular iteration intervals. 
A simulated annealing-type acceptance criterion allows occasional acceptance of worse solutions to escape local minima. 
The ALNS terminates after a fixed number of iterations, returning an initial routing plan
\begin{equation}
\sigma^0 = \{ \pi_w^0 \}_{w \in \mathcal{W}_{\mathrm{idle}}(0)}
\end{equation}
for the known tasks and the workers available at time~0. 
At this stage, the plan is a static solution over $K_0$ and the workers available at time~0. 
Thereafter, routing decisions are updated in an event-driven rolling-horizon manner when new tasks arrive and when workers become idle (Section~\ref{subsubsec:event_update}).

\subsubsection{Event-driven time update and myopic dispatch}
\label{subsubsec:event_update}

During execution, the routing plan is updated in an event-driven manner: time does not advance in fixed increments, but jumps from one decision-relevant event to the next. 
We distinguish three types of events:
\begin{itemize}
	\item \emph{Task arrival}: a new task $i^\ast$ with
	release time $r_{i^\ast}$ becomes known.
	\item \emph{Worker becomes idle}: a worker $w$ becomes available for dispatch at time $t$.
	This occurs either (i) when $t = T_{start}^w$ and the worker starts its shift and is idle at its origin $s_w$, or
	(ii) when the worker completes service at its current task and becomes idle, ready to travel to the next node on its route.
	\item \emph{Horizon termination}: the end of the planning
	horizon, by which all workers must return to their destinations.
\end{itemize}
Let $0 = \mathbb{T}_0 < \mathbb{T}_1 < \mathbb{T}_2 < \dots$ denote the increasing sequence of times at which one or more of these events occur. 
Each $\mathbb{T}_k$ defines a decision epoch at which the routing plan may be revised.

At a generic decision epoch with current time $t = \mathbb{T}_k$, the dispatcher first updates the dynamic state $X(t)$ as defined in Section~\ref{sec:problem}, in particular the set of available tasks $A(t)$, the executed prefixes $R_w^t$ of all worker routes, and the set of idle workers $\mathcal{W}_{\mathrm{idle}}(t)$. 
Workers that are currently travelling to or serving a previously assigned task are treated as non-idle; their routes and immediate next destinations are kept unchanged until they complete the current service.

A static HT-TOPTW subproblem is then formulated on the currently available tasks and idle workers. The node set is
\begin{equation}
V(t) = A(t) \cup V_S(t) \cup V_D(t),
\end{equation}
where $V_S(t) = \{s_w(t)\}_{w \in \mathcal{W}_{\mathrm{idle}}(t)}$ and $V_D(t) = \{d_w\}_{w \in \mathcal{W}_{\mathrm{idle}}(t)}$. 
For each idle worker $w \in \mathcal{W}_{\mathrm{idle}}(t)$, the origin $s_w(t)$ represents its current location at time $t$, and its temporal parameters are updated to
\[
T_{start}^w(t) = t, \qquad T_{end}^w(t) = T_{end}^w.
\]
Tasks $i \in A(t)$ retain their profits, service durations, time windows, and release times. 
Non-idle workers and already executed route segments $\{R_w^t\}$ are treated as fixed commitments and are not modified by the optimization at epoch $t$.

Because workers have heterogeneous origins, destinations, and time budgets, many available tasks are infeasible for all currently idle workers at a given decision epoch. 
Infeasibility arises even under immediate departure. 
Specifically, a worker may be unable to reach the task within its time window while still returning to its destination by the deadline.
Keeping such tasks in the snapshot only enlarges the search space without affecting the set of feasible near-term actions. 
To exploit this structure, we introduce a simple but problem-specific feasibility pre-screening step before invoking the ALNS at each decision epoch.

Let $\mathcal{W}_{\mathrm{idle}}(t)$ denote the set of idle workers and $A(t)$ the set of currently available tasks at time $t$, as defined above. For each task $i \in A(t)$ and each idle worker $w \in \mathcal{W}_{\mathrm{idle}}(t)$, we compute a lower bound on the service start time at $i$,
\begin{equation}
\underline{a}^w_i(t) \;=\;
\max\bigl\{\, t + t_{s_w(t),i},\; r_i,\; b_i \,\bigr\},
\end{equation}
where $s_w(t)$ is the current location of worker $w$ at time $t$, $t_{s_w(t),i}$ is the travel time from $s_w(t)$ to $i$, $r_i$ is the release time of task $i$, and $[b_i,e_i]$ is its time window. 
A necessary condition for worker $w$ to be able to serve task $i$ and still return to its destination by the deadline $T_{end}^w$ is
\begin{equation}
\underline{a}^w_i(t) \leq e_i \qquad \underline{a}^w_i(t) + \tau_i + t_{i,d_w} \;\le\; T_{end}^w,
\end{equation}
where $\tau_i$ is the service duration at $i$ and $t_{i,d_w}$ is the travel time from $i$ to $d_w$. 
We then define the set of idle workers for which task $i$ is potentially feasible at time $t$ as
\begin{equation}
F_i(t) \;=\;
\left\{\, w \in \mathcal{W}_{\mathrm{idle}}(t) \;:\;
\underline{a}^w_i(t) \leq e_i,\; \underline{a}^w_i(t) + \tau_i + t_{i,d_w} \le T_{end}^w \,\right\}.
\end{equation}
If $F_i(t) = \emptyset$, task $i$ cannot be served by any currently idle worker without violating its deadline and is therefore omitted from the snapshot at epoch $t$. 
Tasks that are filtered out in this way remain part of the global task pool and may re-enter future snapshots once additional workers become idle. 
This pre-screening step leaves the myopic structure of the policy unchanged. 
At each epoch, it exploits heterogeneity in worker locations and time budgets to quickly filter globally infeasible tasks. 
This filtering substantially reduces the size of the resulting static subproblems.

The resulting static HT-TOPTW subproblem over $V(t)$ and $\mathcal{W}_{\mathrm{idle}}(t)$ is then approximately solved by the same ALNS as in the initial phase, but with a limited iteration budget appropriate for real-time use. 
The ALNS produces updated routes $\{ \tilde{\pi}_w^t \}_{w \in \mathcal{W}_{\mathrm{idle}}(t)}$ for the idle workers, each starting at $s_w(t)$ at time $t$ and ending at $d_w$ before $T_{end}^w$, and visiting a subset of tasks in $A(t)$.

From these updated routes, the dispatcher extracts near-term actions to be implemented before the next decision epoch. 
A simple myopic dispatch rule is used as a baseline in our experiments. It commits at most one new task to each idle worker. 
Specifically, for each $w \in \mathcal{W}_{\mathrm{idle}}(t)$, if $\tilde{\pi}_w^t = (s_w(t), i_1, i_2, \dots, d_w)$ contains at least one task node, then the dispatcher assigns the first task $i_1$ to worker $w$ as its next destination.
The worker then travels from its current location to $i_1$, waits if necessary to respect $[b_{i_1},e_{i_1}]$ and $r_{i_1}$, serves the task, and subsequently becomes idle again, triggering a new decision epoch. 
If the updated route contains no task nodes for a given worker, that worker either waits at its current location until a new task arrives or returns directly to its destination $d_w$ if doing so is necessary to respect the deadline $T_{end}^w$.

This rolling-horizon ALNS is myopic in the sense that each static subproblem considers only tasks that are currently available in $A(t)$ and ignores the impact of potential future arrivals. 
It nevertheless exploits the full structure of HT-TOPTW at each decision epoch and already provides a strong dynamic baseline. 
In Section~\ref{subsec:scenario_extension}, we build on this myopic scheme by embedding it into a scenario-sampling framework that incorporates sampled future demand and produces joint dispatch decisions for multiple workers based on multi-scenario aggregation.

\subsection{Scenario-based extension and multi-worker joint dispatch}
\label{subsec:scenario_extension}

The rolling-horizon ALNS of Section~\ref{subsec:rh_alns} optimizes only with respect to tasks that are currently available in $A(t)$ and treats future arrivals as unknown. 
As in many dynamic routing settings, such a myopic policy can be short-sighted. First, it may commit workers to attractive tasks far from their current operating area, which reduces their flexibility for future requests. Second, it may consume most of the remaining time budget near the end of the horizon.
To alleviate this limitation, we augment the snapshot-based optimization with a scenario-sampling lookahead mechanism.

The core mechanism is to approximate the impact of future task arrivals by sampling a set of possible future demand realizations at each decision epoch. 
For a given epoch $t$, we generate $S$ independent scenarios of future requests over $(t, H]$ and, in each scenario, augment the current set of real tasks $A(t)$ with a small number of virtual tasks that represent potential future requests. 
Each augmented instance is then solved approximately by the static HT-TOPTW ALNS described in Section~\ref{subsec:rh_alns}. 
From each scenario solution, we extract for every idle worker at most one candidate first real task that is feasible to serve next in that scenario. 
These per-scenario candidate actions are finally aggregated across scenarios by a consensus-type rule to produce a single joint dispatch decision, assigning at most one task to each idle worker and at most one worker to each task.

This procedure yields a scenario-sampling rolling-horizon policy. 
Compared with the myopic ALNS, it biases near-term assignments toward route structures that are repeatedly favored across sampled futures. 
This produces next-task decisions that are more robust to uncertainty in future demand while preserving the original computational structure.
The following subsections describe the scenario generation, virtual-task construction, and multi-scenario aggregation steps in more detail.

\subsubsection{Scenario generation and virtual tasks}
\label{subsubsec:scenario_generation}

At a generic decision epoch with current time $t$ and dynamic state
\begin{equation}
X(t) = (A(t), \allowbreak \{R_w^t, \pi_w^t\}_{w \in \mathcal{W}})
\end{equation}
as defined in
Section~\ref{sec:problem}, we construct $S$ independent demand scenarios. 
Each scenario $s = 1,\dots,S$ is defined by an augmented task set
\begin{equation}
\widetilde{A}^{(s)}(t)
\;=\; A(t) \cup V^{(s)}(t),
\end{equation}
where $V^{(s)}(t)$ is a set of $N^{\mathrm{vir}}$ virtual tasks that act
as stochastic proxies for possible future requests beyond time $t$. 
Virtual tasks are artificial entities used only within the optimization procedure; they are never actually dispatched to workers and are discarded after the dispatch decision at epoch $t$ has been made.

We adopt a lightweight, data-driven scheme for sampling virtual-task attributes. 
Let $\mathcal{L}(t)$ denote the set of coordinates of the currently available real tasks in $A(t)$ and the origin and destination locations of all workers. 
We compute a bounding box enclosing $\mathcal{L}(t)$ and draw the location of each virtual task uniformly at random within this box.
This ensures that virtual tasks populate the same geographical region as the real system. 
For profits and service durations, we sample from empirical ranges inferred from the already known real tasks. 
Specifically, letting $p_i$ and $\tau_i$ denote the profit and service duration of real tasks $i \in A(t)$, we draw the profit $\tilde{p_j}$ of a virtual task $j$ from a uniform distribution on $[\min_{i \in A(t)} p_i,\, \max_{i \in A(t)} p_i]$ and its service duration $\tilde{\tau}_j$ from a uniform distribution on $[\min_{i \in A(t)} \tau_i,\, \max_{i \in A(t)} \tau_i]$. 
This preserves, in a coarse way, the scale of rewards and service times observed in the current instance while avoiding strong parametric assumptions on the future demand.

Time-window parameters for virtual tasks are sampled in two steps.
First, for each virtual task $j$ we draw a time-window lower bound $\tilde{b}_j$ uniformly over the remaining horizon $[t, H]$. 
Second, we sample a time-window width $\Delta$ from the empirical range of widths $e_i - b_i$ observed among real tasks, and set the latest service time to
\begin{equation}
\tilde{e}_j = \min\{\tilde{b}_j + \Delta,\;
H\}.
\end{equation}
The interval $[\tilde{b}_j, \tilde{e}_j]$ is then used directly as the time window of the virtual task in the static HT-TOPTW instance: service cannot start before $\tilde{b}_j$ and must start no later than $\tilde{e}_j$.

In summary, each virtual task $j \in V^{(s)}(t)$ is characterized by a location $\tilde{\ell}_j$, a profit $\tilde{p}_j$, a service duration $\tilde{\tau}_j$, and a time window $[\tilde{b}_j, \tilde{e}_j]$. 
Together with the current real tasks $A(t)$ and the idle workers at time $t$, these virtual tasks define an augmented static HT-TOPTW instance for
scenario $s$. 
In Section~\ref{subsubsec:scenario_solution} we explain how this augmented instance is solved by ALNS and how the solution is used to extract candidate first real tasks for each idle worker.

\subsubsection{Solving scenarios and extracting per-worker candidate actions}
\label{subsubsec:scenario_solution}

For each scenario $s = 1,\dots,S$ constructed in Section~\ref{subsubsec:scenario_generation}, we define an augmented static HT-TOPTW instance at time $t$ with task set 
\begin{equation}
	\widetilde{A}^{(s)}(t) = A(t) \cup V^{(s)}(t),
	\label{eq:aug_taskset}
\end{equation}
and idle workers $\mathcal{W}_{\mathrm{idle}}(t)$. Non-idle workers and already executed route segments $\{R_w^t\}$ are treated as fixed commitments and are not modified.
Formally, the node set is
\begin{equation}
	\widetilde{V}^{(s)}(t)
	= \widetilde{A}^{(s)}(t) \cup V_S(t) \cup V_D(t),
	\label{eq:aug_nodeset}
\end{equation}
with $V_S(t)$ and $V_D(t)$ defined as in Section~\ref{subsec:rh_alns}. Worker-specific start times are $T_{start}^w(t)=t$ for all $w \in \mathcal{W}_{\mathrm{idle}}(t)$, deadlines remain $T_{end}^w$, and travel times $t_{ij}$ are unchanged. 
Profits, service durations, release times and time windows of real tasks in $A(t)$ are kept identical to the original problem data; virtual tasks $j \in V^{(s)}(t)$ use the sampled attributes \[(\tilde{p}_j,\tilde{\tau}_j, [\tilde{b}_j,\tilde{e}_j]).\]

On this augmented instance, we run the same ALNS as in Section~\ref{subsec:rh_alns} to obtain a set of scenario-specific routes
\begin{equation}
	\widetilde{\sigma}^{(s)}(t)
	= \{\, \widetilde{\pi}_w^{(s)}(t) \,\}_{w \in \mathcal{W}_{\mathrm{idle}}(t)}.
	\label{eq:scenario_routes}
\end{equation}

Each route $\widetilde{\pi}_w^{(s)}(t)$ starts at the worker's current location $s_w(t)$ at time $t$, ends at $d_w$ before $T_{end}^w$, and visits a mixture of real tasks from $A(t)$ and virtual tasks from $V^{(s)}(t)$ in some order. 
Because virtual tasks are treated symmetrically to real ones inside ALNS, they can influence the structure of these routes by attracting or repelling the worker in space and time, even though they will never be executed.

Starting from the worker's actual state at time $t$, we simulate travel along this route and identify the first real task $i \in A(t)$ that can be reached and served while satisfying all temporal constraints. 
More precisely, we recompute the service start times along $\widetilde{\pi}_w^{(s)}(t)$ using the timing recursions introduced in Section~\ref{sec:problem}. 
The initial condition is given by the worker’s state at time~$t$. 
We then retain as a candidate the first real task whose time-window and deadline constraints are satisfied.
We ignore all virtual tasks during this extraction step: if $j_g$ corresponds to a virtual task, it is skipped and the simulation continues to $j_{g+1}$. 
The first real task $i$ encountered that satisfies all temporal constraints is recorded as the per-scenario candidate action for worker $w$ in scenario $s$. 
If no such real task exists on $\widetilde{\pi}_w^{(s)}(t)$, then scenario $s$ does not propose a candidate for $w$.

In summary, solving scenario $s$ yields a set of worker-task pairs
\begin{equation}
	C^{(s)}(t) \subseteq
	\mathcal{W}_{\mathrm{idle}}(t) \times A(t),
	\label{eq:scenario_candidates}
\end{equation}
where each $(w,i) \in C^{(s)}(t)$ indicates that, under scenario $s$, sending worker $w$ to task $i$ next is both feasible and recommended by the ALNS route structure. 
Within a given scenario $s$, each idle worker appears in $C^{(s)}(t)$ at most once. 
Across different scenarios, however, the same real task $i$ may be recommended as the first task for different workers.

\begin{figure*}[!tbp]
	\centering
	\subfigure[ALNS solution on augmented task set $\widetilde{A}^{(s)}(t)=\{1,2,3,v\}$.]{
		\label{fig:vtask_a}
		\includegraphics[width=0.32\linewidth]{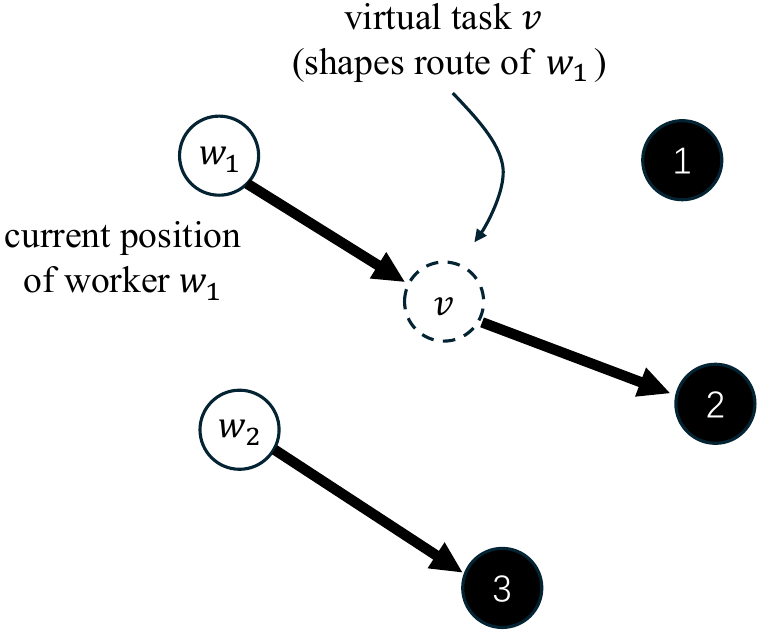}}
	\hspace{2cm}  
	\subfigure[The virtual task $v$ is ignored when extracting $C^{(s)}(t)$.]{
		\label{fig:vtask_b}
		\includegraphics[width=0.24\linewidth]{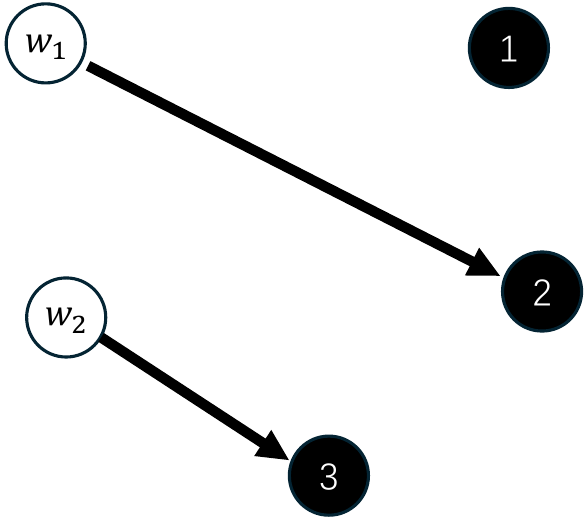}}
	\caption{Example of scenario-specific worker-task extraction.} 
	\label{fig:vtask}
\end{figure*}

\textbf{Example:} For illustration, consider a decision epoch with two idle workers $w_1$ and $w_2$ and three currently available \emph{real} tasks $A(t)=\{1,2,3\}$. 
In a given scenario $s$, we augment this set with a virtual task $v\in V^{(s)}(t)$, so that the ALNS is solved on the augmented task set 
\[\widetilde{A}^{(s)}(t)=A(t)\cup V^{(s)}(t)=\{1,2,3,v\}.\] Suppose that in this scenario, the ALNS solution routes worker $w_1$ from its current position to the virtual task $v$ and then toward real task~2, while worker $w_2$ is routed directly toward real task~3. 
When applying the extraction procedure described above, we ignore $v$ and identify task~2 as the first feasible real task for $w_1$ and task~3 for $w_2$, yielding \[C^{(s)}(t)=\{(w_1,2),(w_2,3)\},\] even though the virtual task $v$ appears first along $w_1$'s route. 
In this way, virtual tasks can shape the routes and thereby influence which real tasks are recommended, but only real tasks are ever retained as candidate actions. 
This example is illustrated in Figure~\ref{fig:vtask}.

\subsubsection{Multi-scenario aggregation and dispatch decision}
\label{subsubsec:aggregation}

The sets $C^{(s)}(t)$ produced by the $S$ scenarios provide up to $S$ recommendations per idle worker, each reflecting one possible realization of future demand. 
To derive a single implementable decision at time $t$, we aggregate these recommendations by counting the frequency with which each worker-task pair appears across scenarios, and then selecting a conflict-free set of pairs with the highest support.

Let
\begin{equation}
	\mathcal{C}(t)
	= \bigcup_{s=1}^S C^{(s)}(t)
	\;\subseteq\;
	\mathcal{W}_{\mathrm{idle}}(t) \times A(t),
	\label{eq:union_candidates}
\end{equation}
be the union of all candidate pairs, and define a frequency map
$f_t : \mathcal{W}_{\mathrm{idle}}(t) \times A(t) \to \{0,1,\dots,S\}$ by
\begin{equation}
	f_t(w,i)
	= \bigl|\{\, s \in \{1,\dots,S\} : (w,i) \in C^{(s)}(t) \,\}\bigr|.
	\label{eq:freq_map}
\end{equation}
Intuitively, $f_t(w,i)$ counts how often sending worker $w$ to task $i$ next is recommended by the scenario solutions. 
Pairs with larger $f_t(w,i)$ are more robust in the sense that they are preferred under many different realizations of future demand.

To construct a dispatch decision, we first discard pairs whose support is too weak. 
Specifically, we retain only those $(w,i) \in \mathcal{C}(t)$ such that
\begin{equation}
	f_t(w,i) \ge \theta_{\min}, \qquad
	\theta_{\min} = \max\{1, \lfloor \alpha S \rfloor\},
	\label{eq:theta_min}
\end{equation}
for some threshold $\alpha \in (0,1)$. 
This ensures that very rarely recommended actions, which may be attributable to sampling noise in individual scenarios, are not selected. 
Among the remaining pairs, we then choose a maximal conflict-free subset by a greedy procedure.

We sort all retained pairs in non-increasing order of their frequency $f_t(w,i)$, breaking ties by a secondary criterion such as larger task profit $p_i$ or shorter travel distance from the worker's current location to task $i$. 
We then traverse this sorted list once and greedily construct the dispatch set $D(t)$, adding a pair $(w,i)$ only if neither worker $w$ nor task $i$ has been selected earlier.
This procedure yields a set of dispatch decisions that satisfies
\begin{equation}
	D(t) \subseteq \mathcal{W}_{\mathrm{idle}}(t) \times A(t),
	\label{eq:dispatch_set}
\end{equation}
in which each worker $w$ and each task $i$ appears at most once.

The set $D(t)$ thus defines a joint dispatch decision: each idle worker $w$ assigned in $D(t)$ is sent to exactly one next task $i$, and each task $i$ is assigned to at most one worker.

\textbf{Example:} Suppose that at a given epoch we have two idle workers $w_1,w_2$, two available tasks $1,2$, and $S=10$ scenarios. 
If the scenario solutions produce counts $f_t(w_1,1)=8$, $f_t(w_1,2)=1$, $f_t(w_2,1)=2$, $f_t(w_2,2)=7$, and $\alpha=0.2$ so that $\theta_{\min}=2$, then the retained pairs are $(w_1,1)$, $(w_2,1)$ and $(w_2,2)$. 
Sorting by frequency yields the order $(w_1,1)$, $(w_2,2)$, $(w_2,1)$. 
The greedy selection then constructs $D(t) = \{(w_1,1),(w_2,2)\}$, assigning each worker to a distinct task that is highly supported across scenarios.

Once $D(t)$ has been determined, the corresponding worker-task assignments are committed in the real system. 
Each assigned worker $w$ departs toward its assigned task $i$ and serves it if the time-window, release-time, and deadline constraints permit. 
The worker then becomes idle again, triggering a new decision epoch.
Workers in $\mathcal{W}_{\mathrm{idle}}(t)$ that do not appear in $D(t)$ either wait at their current locations until new tasks arrive or return to their destinations if necessary to satisfy their deadlines. 
In particular, if a worker's candidate assignments receive little or conflicting support across scenarios, the worker is likely to be omitted from $D(t)$, so the policy effectively chooses to wait for additional information rather than committing to a poorly supported task. 
In this way, the scenario-based extension preserves the rolling-horizon structure of the  myopic ALNS, but biases near-term decisions toward worker-task assignments that are repeatedly supported across a diverse set of sampled futures.

\subsection{Overall algorithm and computational aspects}
\label{subsec:overall}

The ingredients introduced in Sections~\ref{subsec:rh_alns} and~\ref{subsec:scenario_extension} together define a complete online control policy for DTOP-SC. 
At a high level, the policy maintains, for each worker, a route consisting of an executed prefix and an unexecuted suffix. 
It monitors the system state through event-driven decision epochs. 
At each epoch, it constructs and solves a family of static HT-TOPTW subproblems to decide which next tasks to assign to idle workers.
The myopic rolling-horizon ALNS of Section~\ref{subsec:rh_alns} is recovered as a special case when no virtual tasks are used and only a single scenario is considered.

Algorithm~\ref{alg:rh_alns} summarizes the Scen-RH-ALNS.  
The algorithm starts from the initial routing plan obtained by solving the static HT-TOPTW instance for the tasks known at time zero. 
During execution, time advances from one event to the next: whenever a new task arrives or a worker completes service and becomes idle, the current time is updated, the dynamic state is recomputed, and a new decision epoch is triggered. 
At each epoch, $S$ augmented static instances are built by adding virtual tasks to the set of currently available real tasks, one instance per scenario. 
Each static instance is solved approximately by ALNS, yielding a set of scenario-specific routes for the idle workers. 
From these routes, at most one candidate first real task is extracted per worker and scenario. 
Finally, the multi-scenario aggregation procedure of Section~\ref{subsubsec:aggregation} is applied to obtain a conflict-free set of worker-task assignments, which are committed in the real system until the next decision epoch.

\begin{algorithm}[t]
	\caption{Scenario-sampling RH-ALNS policy for DTOP-SC}
	\label{alg:rh_alns}
	\begin{algorithmic}[1]
		\State Initialize $t \gets 0$; construct initial task set $K_0$ and static HT-TOPTW instance
		\State Solve initial static instance by ALNS to obtain routes $\sigma^0 = \{\pi_w^0\}_{w \in \mathcal{W}_{\mathrm{idle}}(0)}$
		\While{$t < H$ and (unserved tasks remain or some worker has not reached $d_w$)}
		\State Advance $t$ to the time of the next event (task arrival, worker becomes idle, or $H$)
		\State Update dynamic state: available tasks $A(t)$, executed prefixes $R_w^t$, idle workers $\mathcal{W}_{\mathrm{idle}}(t)$
		\If{$t = H$ or $\mathcal{W}_{\mathrm{idle}}(t) = \emptyset$ or $A(t) = \emptyset$}
		\State Continue execution of already committed routes; \textbf{continue}
		\EndIf
		\For{$s = 1,\dots,S$}
		\State Generate virtual tasks $V^{(s)}(t)$; form augmented task set $\widetilde{A}^{(s)}(t)$ as in~\eqref{eq:aug_taskset}
		\State Build augmented static HT-TOPTW instance for scenario $s$ at time $t$ with node set $\widetilde{V}^{(s)}(t)$ in~\eqref{eq:aug_nodeset}
		\State Run ALNS on scenario $s$ to obtain routes $\widetilde{\sigma}^{(s)}(t)$ for idle workers as in~\eqref{eq:scenario_routes}
		\State Extract candidate worker-task pairs $C^{(s)}(t)$ (first feasible real task per worker) as in~\eqref{eq:scenario_candidates}
		\EndFor
		\State Aggregate $\{C^{(s)}(t)\}_{s=1}^S$ using~\eqref{eq:union_candidates}--\eqref{eq:theta_min} to compute dispatch decision $D(t)$ in~\eqref{eq:dispatch_set}
		\State Commit assignments in $D(t)$; update planned routes and executed prefixes accordingly
		\EndWhile
	\end{algorithmic}
\end{algorithm}

From a computational standpoint, the main cost of the policy at each decision epoch comes from solving the $S$ augmented static HT-TOPTW instances by ALNS. 
Let $I_{\text{ALNS}}$ denote the number of ALNS iterations per scenario, and let $m_t = |\mathcal{W}_{\mathrm{idle}}(t)|$ and $n_t = |A(t)|$ be the numbers of idle workers and available tasks in the snapshot at epoch $t$. 
In our implementation, each iteration performs a constant number of destroy/repair operations and local search moves whose cost is polynomial in the total route length, leading to an overall per-epoch effort on the order of
\[
O\bigl(S \, I_{\text{ALNS}} \, (m_t + n_t)^3\bigr)
\]
in the worst case, dominated by regret-$k$ insertion and local-search moves. 
Because the static instances at a given epoch are independent across scenarios, this cost is highly parallelizable: each scenario can be allocated to a separate processing unit, and the overall latency is then dominated by the slowest scenario. 
However, this structure is exploited by distributing scenarios across multiple cores and solving them in parallel, so that the wall-clock time per epoch grows sublinearly with~$S$.

Several design choices help keep the per-epoch computation compatible with real-time operation. 
First, the static instances are defined only on the currently available tasks $A(t)$ and idle workers $\mathcal{W}_{\mathrm{idle}}(t)$, which are typically much smaller than the full task set and workforce. 
Second, the number of scenarios $S$ and the number of virtual tasks per scenario $N^{\mathrm{vir}}$ are treated as algorithmic parameters: larger values enhance the fidelity of the lookahead but increase computation, whereas smaller values lead to faster but more myopic decisions. 
In the numerical study, we report the empirical trade-off between these parameters and both solution quality and running time.

\section{Computational experiments}
\label{sec:experiments}

This section evaluates the proposed Scen-RH-ALNS policy on both map-based DTOP-SC instances and on the standard DTOP benchmark. 
We first describe the instance families, implementation details, and evaluation metrics, and then report numerical results and their interpretation.

\subsection{Experimental setup}
\label{subsec:setup}

This subsection describes two complementary sets of instances, namely synthetic map-based DTOP-SC instances reflecting our target application and the DTOP benchmark of \citet{kirac2025} as an external reference. 
Additionally, we provide details on the Scen-RH-ALNS implementation and shared performance metrics.

\subsubsection{DTOP benchmark and specialization}
\label{subsubsec:dtop_data}

The first family of instances is based on the DTOP benchmark introduced by \citet{kirac2025}. 
In this benchmark, a homogeneous fleet of $2$--$4$ vehicles operates from a single depot over a fixed planning horizon. 
Customer requests are revealed dynamically according to three degrees of dynamism (weak, medium, and high). 
These settings specify the fraction of requests that become known progressively over time.
The seven problem sets p1--p7 contain between 30 and 100 requests (from p1 with 30 requests up to p7 with 100), yielding 1161 instances in total. Section~\ref{subsubsec:dtop_mapping} describes how they are embedded into the DTOP-SC framework. 
Further details on the instance generation procedure can be found in \citet{kirac2025}.

\subsubsection{Map-based DTOP-SC instances}
\label{subsubsec:data}

The second family of instances consists of synthetic DTOP-SC scenarios generated from real-world urban coordinates. 
These instances are designed to capture key structural features of multi-worker field-service systems, including heterogeneous worker trajectories, tight time windows, and dynamically released tasks.

Node coordinates are sampled from a realistic urban map and projected onto a two-dimensional plane. 
Symmetric travel times are set equal to the Euclidean distances between node coordinates.
Worker origins and destinations are drawn independently from the node set subject to a minimum separation of $40\%$ of the map diameter. 
For each worker $w$, we compute the direct origin-to-destination travel time and multiply it by a buffer factor sampled uniformly from $[1.3, 2.5]$ to obtain an individual time budget $T_{\text{avail}}^w$ within a common planning horizon $H = 180$. 
The worker’s start time is drawn uniformly from $[0, \max\{0, H - T_{\text{avail}}^w\}]$ and its deadline is set to the start time plus $T_{\text{avail}}^w$. 
This construction yields a heterogeneous team with different origin/destination locations and effective working-time windows.

Tasks are generated conditionally on the worker ensemble by rejection sampling, so that each accepted task is individually feasible for at least one worker. 
For each candidate task, we sample a location from the node set, a profit from a prescribed profit range, and a service duration from a given duration range. 
We then draw an appearance time $\text{apr}$ from a feasible working interval of a randomly selected worker, ensuring that a fixed slack remains before that worker’s deadline.
Given $\text{apr}$, we construct a service time window $[rdy, ldt]$ by sampling $rdy$ uniformly from a feasible interval starting at $\text{apr}$ plus a small minimum-travel buffer, and a window width from the configured range, truncating $ldt$ at the horizon $H$ if necessary.

The candidate task is accepted only if it is time-feasible for at least one worker $w$. 
Time feasibility means that $w$ can travel from its origin to the task, possibly wait, and start service within $[rdy, ldt]$. 
After service, $w$ must still be able to reach its destination before its individual deadline.
This procedure produces tasks that are heterogeneous in spatial location and temporal feasibility and that induce worker-specific local orienteering subproblems. 
All generated tasks are finally pooled and randomly shuffled to remove ordering effects.

For numerical stability and to reuse a single set of heuristic parameters across heterogeneous instances, we rescale task profits in each map-based instance as $\tilde p_i = p_i / p_{\mathrm{scl}}$, where $p_{\mathrm{scl}}$ is the upper bound of the profit range used by the generator for that instance family.
Throughout this section, profit ranges such as $[10,50]$ refer to the pre-scaling values used by the generator; all computations and reported profits use the scaled values.
Accordingly, all map-based profits reported in Section~\ref{subsec:roadnet_results} are in this scaled unit.
This positive linear scaling does not affect relative comparisons, and the factor cancels out in the gap metric $\mathrm{Gap}_{MIP}$.

Across this map-based setting we define a \emph{base} configuration and several one-factor perturbations. 
In the base configuration, we consider $M=10$ workers and $N=100$ tasks; task service durations $\tau_i$ are drawn uniformly from $[1,3]$ time units, time-window widths are drawn from $[10.0,20.0]$, and task profits $p_i$ are drawn from $[10,50]$. 
To study sensitivity to service times, we create a \emph{Short} duration family with $\tau_i \in [0,2]$ and a \emph{Long} duration family with $\tau_i \in [2,6]$, keeping all other parameters at their base values.
To examine the effect of temporal flexibility, we define \emph{Tight} and \emph{Loose} time-window families, in which time-window widths are drawn from $[5.0,15.0]$ and $[15.0,30.0]$, respectively, while $M$, $N$, $\tau_i$ and $p_i$ follow the base configuration. 
Profit heterogeneity is varied through \emph{Narrow} and \emph{Wide} profit families, where profits are drawn from $[10,20]$ and $[10,100]$, respectively, with all other parameters as in the base setting.
Finally, to assess scalability, we construct a \emph{Scale} family in which the number of workers and tasks grow proportionally, using $(M,N) \in \{(5,50), (7,70), (9,90), (11,110), (13,130), (15,150)\}$, and keeping $\tau_i \in [1,3]$, time-window widths in $[10.0,20.0]$, and profits in $[10,50]$. 
For each configuration described above, we generate ten independent instances, which are used to assess the stability and scalability of Scen-RH-ALNS on highly heterogeneous DTOP-SC settings.

\begin{figure}[!tbp]
	\centering\includegraphics[width=0.5\linewidth]
	{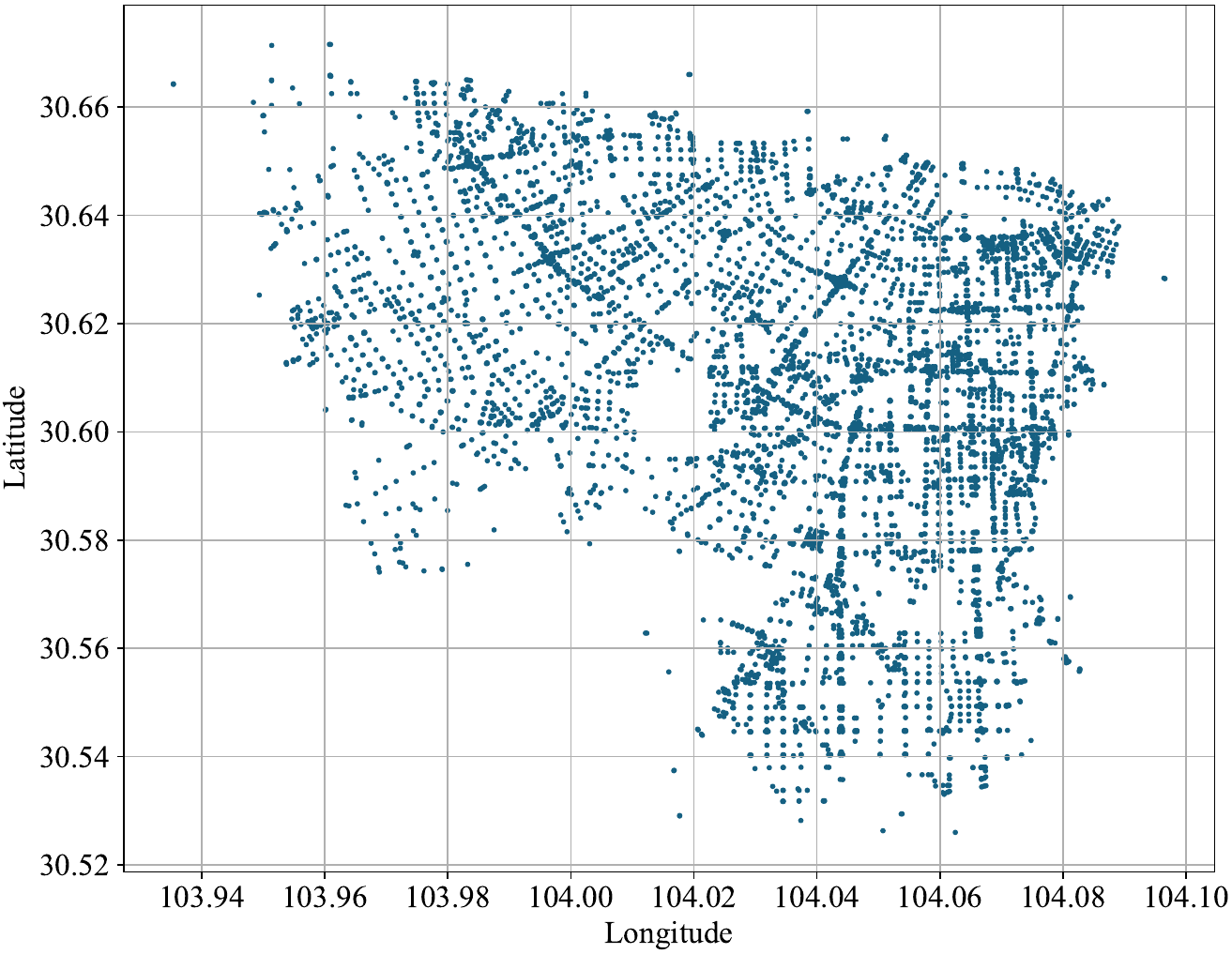}
	\caption{Spatial distribution of sampled road-network coordinates used to generate the synthetic instances.}
	\label{fig:roadnet_nodes}
	\end{figure}

Figure~\ref{fig:roadnet_nodes} visualizes the spatial distribution of the sampled road-network coordinates used in our generator. 
The point cloud exhibits a clearly non-uniform spatial structure characteristic of urban areas, with several dense corridors and sparser peripheral regions, which is typical of real city layouts.

\subsubsection{Implementation details and parameter settings}
\label{subsubsec:impl}

All experiments are implemented in Python using a common code base for both instance families. 
The rolling-horizon ALNS described in Section~\ref{sec:method} is used without instance-specific tuning: the same destroy/repair operators, local search moves, and acceptance scheme are employed for all static HT-TOPTW snapshots and all scenario-augmented instances.

For Scen-RH-ALNS, we fix the number of demand scenarios at $S = 15$, the number of virtual tasks per scenario at $N^{\mathrm{vir}} = 5$, and the threshold $\alpha = 0.2$. 
Within each scenario, the corresponding static augmented instance is solved by ALNS with an iteration limit $I_{\text{ALNS}}^{\text{scen}} = 100$. 
At the beginning of each horizon, the initial static HT-TOPTW instance is solved with a larger limit $I_{\text{ALNS}}^{\text{init}}$ to construct a reasonably good starting solution, consistent with Section~\ref{sec:method}. 
This independence allows us to evaluate the $S$ augmented instances in parallel using a multi-process executor, assigning each scenario to a separate CPU core whenever possible.
This reduces the wall-clock time per decision epoch and makes the scenario-sampling extension compatible with real-time operation.

Unless otherwise stated, all runtimes reported for our algorithms are elapsed wall-clock times, measured on a workstation equipped with an Intel Core i5 CPU at 3.5\,GHz.
For DTOP benchmark instances, we report times only for our rolling-horizon ALNS variants. 
The times for MPAc and MPAd are taken from the tables reported in \citet{kirac2025}, which were obtained on a different hardware platform. As a result, these times are not directly comparable in an absolute sense and are used only as indicative references to illustrate orders of magnitude. In addition to Scen-RH-ALNS, we also consider a myopic rolling-horizon variant (Myopic-RH) that uses the same HT-TOPTW snapshot model and ALNS operators but sets the number of scenarios to one and the number of virtual tasks per snapshot to zero.

\subsubsection{Evaluation metrics}
\label{subsubsec:metrics}

We employ related but distinct performance metrics on the two instance families.

\paragraph{DTOP benchmark instances}
On the DTOP benchmark, our primary comparison is in terms of total collected profit. For MPAc and MPAd, we use the published profit values from \citet{kirac2025}. Average times are reported for all methods, but because MPAc/MPAd were executed on different hardware, their times are interpreted only at the level of orders of magnitude relative to our implementations. To complement average profit, we also consider instance-wise best counts, defined as the number of instances on which a given method attains at least as much profit as its competitors, and optimality gaps with respect to the CP-based offline solutions. For each instance with offline reference value $Z_{\text{CP}}$, the percentage optimality gap of a method with objective value $Z$ is computed as
\[
\text{Gap}_{\text{CP}} = 100\% \times \frac{Z_{\text{CP}} - Z}{Z_{\text{CP}}}.
\]
Taken together, these metrics allow us to assess average performance, relative optimality, and the distribution of outcomes across individual DTOP instances.

\paragraph{Map-based DTOP-SC instances}
On the map-based DTOP-SC instances, we evaluate performance relative to an offline, time-limited MIP reference.
For each instance, we construct the corresponding static HT-TOPTW model by revealing all tasks at time $t=0$ and ignoring dynamic information, and solve this MIP with Gurobi under a time limit of 600 seconds.
We record the best feasible incumbent value $Z_{\text{MIP}}$ returned within the time limit.
Because HT-TOPTW is computationally challenging and the instances are heterogeneous, the MIP may (i) return a weak incumbent within 600 seconds or (ii) fail to find any feasible solution on the largest instances.
Therefore, $Z_{\text{MIP}}$ is interpreted as an offline reference under a fixed computational budget (i.e., a strong heuristic baseline when it succeeds).

Given the profit $Z_{\text{S-ALNS}}$ obtained by Scen-RH-ALNS, we define the percentage gap (when $Z_{\text{MIP}} > 0$) as:
\[
\text{Gap}_{\text{MIP}} = 100\% \times \frac{Z_{\text{MIP}} - Z_{\text{S-ALNS}}}{Z_{\text{MIP}}}.
\]
Positive values indicate that the time-limited MIP incumbent is better than the online policy, whereas negative values indicate that Scen-RH-ALNS outperforms the MIP incumbent within the same 600s budget.
Importantly, when $Z_{\text{MIP}}$ is very small, the ratio above becomes numerically unstable and can yield large-magnitude gaps; such cases primarily reflect difficulty in obtaining a strong MIP incumbent within the time limit.
For instances where $Z_{\text{MIP}}=0$ within 600 seconds, we mark $\text{Gap}_{\text{MIP}}$ as `-'.
For each instance family, we summarize performance by reporting the mean and standard deviation of profits, the average $\text{Gap}_{\text{MIP}}$, and average times.
These metrics provide a consistent basis for evaluating stability, scalability, and the distance to offline references under a fixed computational budget.

\subsection{Validation on DTOP benchmark instances}
\label{subsec:dtop_results}

We first evaluate Scen-RH-ALNS on the standard DTOP benchmark of \citet{kirac2025} to validate the behavior of the proposed online control framework. We then evaluate the same framework on our map-based DTOP-SC instance families in Section~\ref{subsec:roadnet_results}.

\subsubsection{Mapping DTOP-SC to the DTOP setting}
\label{subsubsec:dtop_mapping}

To apply Scen-RH-ALNS to the DTOP benchmark instances of \citet{kirac2025},  we use the homogeneous single-depot specialization of DTOP-SC described in Section~\ref{subsec:ht_toptw_mip}. In all instances, workers  are identical vehicles that share a common origin and destination depot, have  the same start time and deadline, and use Euclidean travel times as in the  original benchmark. For static customers we set $r_i = 0$, whereas for dynamic  customers we follow \citet{kirac2025} and set the release time equal to the start of the time window and the latest service start time equal to the end of the  planning horizon. Under this configuration, each DTOP instance is directly interpreted as a DTOP-SC instance on which Scen-RH-ALNS can be applied.

For performance evaluation, we use the CP-based offline solutions and the multiple plan approaches MPAc and MPAd reported by \citet{kirac2025} as external references, relying on their published profit values and indicative times when comparing routing quality and computational effort. The next subsections report profit, runtime, and optimality-gap comparisons under this common experimental protocol.

\subsubsection{Profit and runtime comparison with MPAc/MPAd}
\label{subsubsec:dtop_profit}

\begin{table}[!t]
	\centering
	\caption{Average profit and times on high-dynamism DTOP instances.}
	\label{tab:dtop_summary}
	\scriptsize
	\setlength{\tabcolsep}{2pt}
	\renewcommand{\arraystretch}{1}
	\begin{tabular}{rlcccc}
		\toprule
		Dynamism & Method & Mean Profit & Mean time (s) & $\Delta$Profit vs MPAd & $\Delta$Profit vs MPAc \\
		\midrule
		\multicolumn{1}{l}{High} & MPAd   & 74.92 & 197.74 &    -   & - \\
		                         & MPAc   & 76.81 & 191.96 &    -   & - \\
		                         & S-ALNS & 74.33 & 0.14   & 0.79\% & 3.23\% \\		
		\bottomrule
	\end{tabular}%
\end{table}

Table~\ref{tab:dtop_summary} summarizes average profits and times under high dynamism for Scen-RH-ALNS and the two multiple plan approaches. In this regime, the mean profit achieved by Scen-RH-ALNS (74.33) is close to that of MPAd (74.92) and MPAc (76.81). Interpreting the percentage differences in the last two columns as gaps relative to MPAd and MPAc, Scen-RH-ALNS remains within about $0.8\%$ of MPAd and $3.2\%$ of MPAc on average. At the same time, the reported average time of Scen-RH-ALNS is around $0.14$ seconds per instance, whereas the times of MPAc and MPAd are approximately 192 and 198 seconds, respectively. Although these times were obtained on different hardware platforms and are therefore not strictly comparable in absolute terms, they are indicative of a difference of roughly two to three orders of magnitude in computational effort.

\begin{table}[htbp]
	\centering
	\caption{Instance-wise best counts on the DTOP benchmark.}
	\label{tab:bestcount}
	\scriptsize
	\setlength{\tabcolsep}{6pt}
	\renewcommand{\arraystretch}{1}
	\begin{tabular}{ccccc}
		\toprule
		\multirow{2}[4]{*}{Degree of Dynamism} & \multirow{2}[4]{*}{\#INS} & \multicolumn{3}{c}{\#BEST} \\
		\cmidrule{3-5}
		&       & S-ALNS & MPAc  & MPAd \\
		\midrule
		Highly     & 387   & 316   & 328   & 356 \\
		Moderately & 387   & 236   & 279   & 306 \\
		Weakly     & 387   & 156   & 283   & 313 \\
		All        & 1161  & 708   & 890   & 975 \\
		\bottomrule
	\end{tabular}%
\end{table}

Table~\ref{tab:bestcount} provides an instance-wise perspective by reporting, for each degree of dynamism, the number of instances on which each method obtains at least as much profit as its competitors. Column ``\#INS'' gives the number of instances in each dynamism group, and the ``\#BEST'' columns report, for each method, how many instances it matches or exceeds the best profit among the three methods. As before, S-ALNS denotes the proposed Scen-RH-ALNS policy.

Under high dynamism, Scen-RH-ALNS matches or exceeds the profit of its competitors on 316 out of 387 instances, while MPAc and MPAd do so on 328 and 356 instances, respectively. Aggregated over all 1161 instances, the corresponding counts are 708 for Scen-RH-ALNS, 890 for MPAc, and 975 for MPAd. These figures indicate that, although the multiple plan approaches attain the highest profit more frequently overall, Scen-RH-ALNS achieves competitive routing quality on a substantial fraction of instances. 

Moreover, when compared against the online solutions reported by \citet{kirac2025}, Scen-RH-ALNS attains strictly higher profit than both MPAc and MPAd on 27, 58, and 30 instances in the high-, moderate-, and low-dynamism settings, respectively, thereby establishing new best-known online results (BKs) on these cases.

Taken together, the average and instance-wise comparisons suggest that Scen-RH-ALNS delivers profits that are close to those of MPAc and MPAd on the DTOP benchmark, particularly under high dynamism, while requiring substantially less computational time. This trade-off between solution quality and runtime reflects the design of Scen-RH-ALNS, which performs relatively shallow but frequent reoptimization at each decision epoch.

\subsubsection{Optimality gap analysis vs.\ CP-based offline solutions}
\label{subsubsec:dtop_gaps}

\begin{table}[htbp]
	\centering
	\caption{Average optimality gaps vs.\ CP-based offline solutions on high-dynamism DTOP instances.}
	\label{tab:dtopgap}
	\scriptsize
	\setlength{\tabcolsep}{6pt}
	\renewcommand{\arraystretch}{1}
	\begin{tabular}{lccc}
		\toprule
		Problem Set & S-ALNS & MPAd  & MPAc \\
		\midrule
		p1    & 55.17\% & 55.56\% &  55.56\% \\
		p2    & 70.00\% & 70.00\% &  70.00\% \\
		p3    & 52.35\% & 52.50\% &  52.50\% \\
		p4    & 47.08\% & 47.22\% &  44.75\% \\
		p5    & 47.04\% & 47.49\% &  46.13\% \\
		p6    & 42.88\% & 43.45\% &  43.45\% \\
		p7    & 45.66\% & 42.46\% &  40.95\% \\
		\midrule
		Mean & 51.45\% & 51.24\% &  50.48\% \\
		SD   &  9.17\% &  9.49\% &  10.02\% \\
		\bottomrule
	\end{tabular}%
\end{table}

To relate online performance to offline benchmarks, Table~\ref{tab:dtopgap} reports average percentage optimality gaps with respect to the CP-based offline solutions of \citet{kirac2025} for high dynamism, broken down by problem set. For each problem set $p1$--$p7$, the table lists the mean gap over instances for Scen-RH-ALNS (S-ALNS), MPAd, and MPAc; the last two rows report the mean and standard deviation of these gaps across the seven problem sets.

Across all three methods, the gaps to the offline CP reference are large, with average values around 50\%. This reflects the intrinsic difficulty of the DTOP benchmark: even strong online policies remain far from the offline upper bound when decisions must be taken in real time. The mean row of Table~\ref{tab:dtopgap} shows that the average gap of Scen-RH-ALNS (51.45\%) is within about one percentage point of MPAd (51.24\%) and MPAc (50.48\%), indicating that all three policies are of similar quality when measured against the CP-based offline benchmark.

At the level of individual problem sets, Scen-RH-ALNS sometimes attains slightly smaller gaps than the multiple plan approaches and sometimes slightly larger ones. For example, Scen-RH-ALNS has the smallest or equal gap on p1, p2, p3, p5, and p6, whereas MPAc or MPAd are marginally closer to the offline reference on p4 and p7. In all cases, however, the differences between methods on a given problem set are only a few percentage points, which is small compared with the overall magnitude of the gaps. These patterns reinforce the view that Scen-RH-ALNS, MPAc, and MPAd achieve comparable levels of offline optimality on the DTOP benchmark.

Overall, the gap analysis suggests that Scen-RH-ALNS offers online performance that is broadly in line with that of state-of-the-art multiple plan approaches when measured relative to CP-based offline solutions, while differing primarily in computational footprint and algorithmic structure.

\subsubsection{Discussion}
\label{subsubsec:dtop_discussion}

The DTOP experiments provide an external check on the generality of the proposed control framework. Without any problem-specific tuning beyond the specialization described in Section~\ref{subsubsec:dtop_mapping}, Scen-RH-ALNS attains average profits that are within a few percentage points of the multiple plan approaches MPAc and MPAd, both in aggregate and when measured relative to CP-based offline solutions. In particular, under high dynamism its mean profit differs from that of MPAd and MPAc by less than 1\% and 3\%, respectively, while the corresponding optimality gaps remain at essentially the same level as those of the state-of-the-art methods.

At the same time, the indicative times highlight a substantial difference in computational footprint: Scen-RH-ALNS requires on the order of $10^{-1}$ seconds per instance on the tested platform, whereas the reported running times for MPAc and MPAd are on the order of $10^{2}$ seconds. Although these times were obtained on different hardware and cannot be compared in absolute terms, they point to a consistent advantage in terms of responsiveness and ease of deployment.

Overall, the DTOP benchmark results indicate that Scen-RH-ALNS delivers competitive routing quality on a standard dynamic orienteering benchmark, including in highly dynamic regimes, while requiring substantially less computation time than multiple plan approaches. This supports the applicability of the proposed rolling-horizon ALNS framework beyond its original heterogeneous DTOP-SC setting.

\subsection{Results on map-based DTOP-SC instances}
\label{subsec:roadnet_results}

We evaluate the proposed Scen-RH-ALNS policy on the map-based DTOP-SC instances introduced in Section~\ref{subsubsec:data}. These experiments are intended to assess the consistency of performance across different instance families, the scalability with respect to the numbers of workers and tasks, and the incremental value of scenario-based lookahead relative to a myopic rolling-horizon baseline. Summary statistics are reported in this section, while complete per-instance numerical results are provided in \ref{append2} and \ref{append3}.

\subsubsection{Stability across instance families}
\label{subsubsec:roadnet_robustness}

\begin{table}[!tbp]
	\centering
	\caption{Aggregate performance of Scen-RH-ALNS (S-ALNS) on map-based instances with varying instance characteristics.}
	\label{tab:robust_dist}
	\scriptsize
	\setlength{\tabcolsep}{3pt} 
	\renewcommand{\arraystretch}{1}
	\begin{tabular}{m{1.4cm}<{\centering} m{1.2cm}<{\centering} m{1.8cm}<{\centering} m{2.0cm}<{\centering} m{1.8cm}<{\centering}  m{1.5cm}<{\centering} m{1.5cm}<{\centering}}
		\toprule
		Family & Dataset& Avg. Z$_{MIP}$ & Avg. $Z_{\text{S-ALNS}}$ & Gap$_{MIP}$& SD& Mean Time (s)\\
		\midrule
		& Base  & 42.94 & 40.56  & 5.55\% & 4.82  & 19.06  \\
		\midrule
		\multicolumn{1}{c}{\multirow{2}[2]{*}{Duration}} & Short & 53.61  & 54.30  & -1.28\% & 5.61  & 101.48  \\
		& Long  & 28.43  & 27.22  & 4.24\% & 3.30  & 4.65  \\
		\midrule
		\multicolumn{1}{c}{\multirow{2}[2]{*}{Time Windows}} & Tight & 42.74  & 41.58  & 2.71\% & 4.00  & 17.64  \\
		& Loose & 38.58  & 38.28  & 0.77\% & 3.31  & 24.08  \\
		\midrule
		\multicolumn{1}{c}{\multirow{2}[2]{*}{Profit}} & Narrow & 48.87  & 48.08  & 1.63\% & 5.26  & 21.42  \\
		& Wide  & 37.49  & 35.68  & 4.84\% & 3.43  & 14.30  \\
		\bottomrule
	\end{tabular}%
\end{table}

Table~\ref{tab:robust_dist} reports aggregate performance metrics for Scen-RH-ALNS across the different map-based instance families defined in Section~\ref{subsubsec:data}. 
For each configuration, we report the average MIP incumbent value $Z_{MIP}$, the average profit $Z_{\text{S-ALNS}}$ obtained by Scen-RH-ALNS, the percentage gaps $\text{Gap}_{MIP}$ as defined in Section~\ref{subsubsec:metrics}, the sample standard deviation (SD) of $Z_{\text{S-ALNS}}$ across instances, and the average time in seconds.

Across all families, the average gaps with respect to $Z_{MIP}$ remain moderate. In the Base configuration, Scen-RH-ALNS achieves an average profit that is within about $5.6\%$ of the MIP incumbent. Variations in service duration, time-window width, and profit range lead to average $\text{Gap}_{MIP}$ values typically between $0\%$ and $5\%$, with the Short-duration setting even exhibiting a slightly negative gap, indicating that Scen-RH-ALNS can occasionally outperform the MIP incumbent within the allotted time. These figures suggest that, on the map-based DTOP-SC instances considered, the online policy remains close to the best offline solutions available under a fixed time limit.

Decision times remain moderate across families and reflect the relative difficulty of the instances. The Short-duration setting, where many short tasks can potentially be sequenced within the workers' time budgets, leads to higher computational effort, whereas the Long-duration configuration is easier to solve and requires only a few seconds on average. Overall, these results indicate that Scen-RH-ALNS maintains near-offline performance while adapting gracefully to changes in service durations, time-window widths, and profit distributions.

\subsubsection{Scalability with the number of workers and tasks}
\label{subsubsec:scalability}

\begin{table}[!tbp]
	\centering
	\caption{Scalability of Scen-RH-ALNS (S-ALNS) on map-based instances.}
	\label{tab:scalability}
	\scriptsize
	\setlength{\tabcolsep}{4pt} 
	\renewcommand{\arraystretch}{1}
	\begin{tabular}{m{1.2cm}<{\centering} m{1.8cm}<{\centering} m{2.0cm}<{\centering} m{1.8cm}<{\centering} m{1.5cm}<{\centering} m{1.5cm}<{\centering}}
		\toprule
		Dataset & Avg.  Z$_{MIP}$& Avg.  $Z_{\text{S-ALNS}}$ & Gap$_{MIP}$& SD& Mean Time (s)\\
		\midrule
		5v/50t   & 23.34  & 22.14  & 5.15\%   & 1.86  & 4.05   \\
		7v/70t   & 30.26  & 28.45  & 5.97\%   & 2.92  & 8.97   \\
		9v/90t   & 38.97  & 36.65  & 5.94\%   & 3.99  & 13.92  \\
		10v/100t & 42.94  & 40.56  & 5.55\%   & 4.82  & 19.06  \\
		11v/110t & 45.80  & 45.02  & 1.70\%   & 4.66  & 34.45  \\
		13v/130t & 48.82  & 51.36  & -5.20\%  & 4.15  & 36.43  \\
		15v/150t & 51.26  & 63.38  & -23.64\% & 5.70  & 60.04  \\
		\bottomrule
	\end{tabular}%
\end{table}

\begin{figure*}[!tbp]
	\centering
	\includegraphics[width=0.7\linewidth]{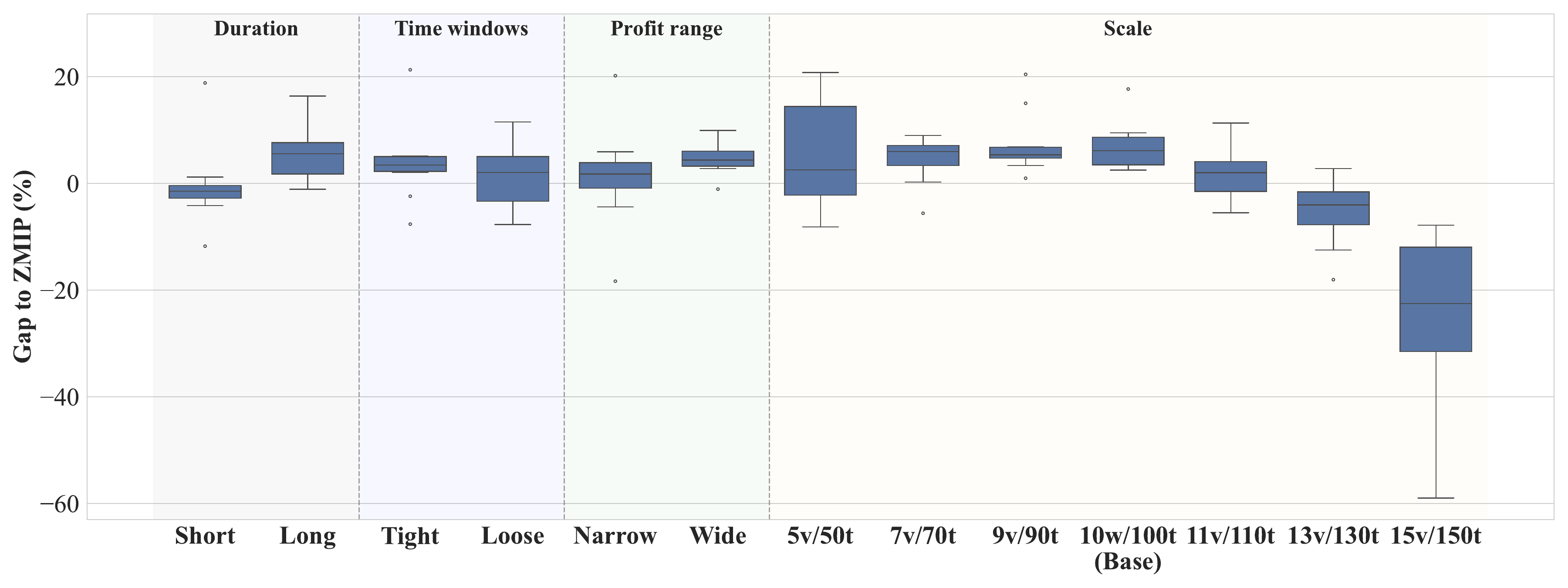}
	\caption{Per-instance distributions of $\text{Gap}_{MIP}$ for map-based stability experiments.}
	\label{fig:roadnet_gaps}
\end{figure*}

Table~\ref{tab:scalability} summarizes the scalability of Scen-RH-ALNS when jointly increasing the number of workers and tasks. For each scale configuration, we report the same set of metrics as in Table~\ref{tab:robust_dist}.

The average gaps with respect to $Z_{MIP}$ remain within a few percentage points at small and medium scales and become negative for the largest instances. For example, in the 5-vehicle/50-task configuration, the average $\text{Gap}_{MIP}$ is about $5.1\%$, while in the largest 15-vehicle/150-task configuration, the average $\text{Gap}_{MIP}$ decreases to approximately $-23.6\%$. This pattern reflects the increasing difficulty of the static HT-TOPTW problems for the MIP solver as the problem size grows, whereas Scen-RH-ALNS continues to produce high-quality solutions within its prescribed iteration limits.

Average times grow smoothly with scale, from roughly 4 seconds for the smallest configuration to about 60 seconds for the largest one. This increase is consistent with the larger HT-TOPTW snapshots and the higher number of decision epochs in larger instances. Taken together, these results indicate that Scen-RH-ALNS scales to realistically sized heterogeneous DTOP-SC instances while maintaining a relatively stable distance to offline benchmarks.

\begin{figure}[!tbp]
	\centering
	\includegraphics[width=0.5\linewidth]{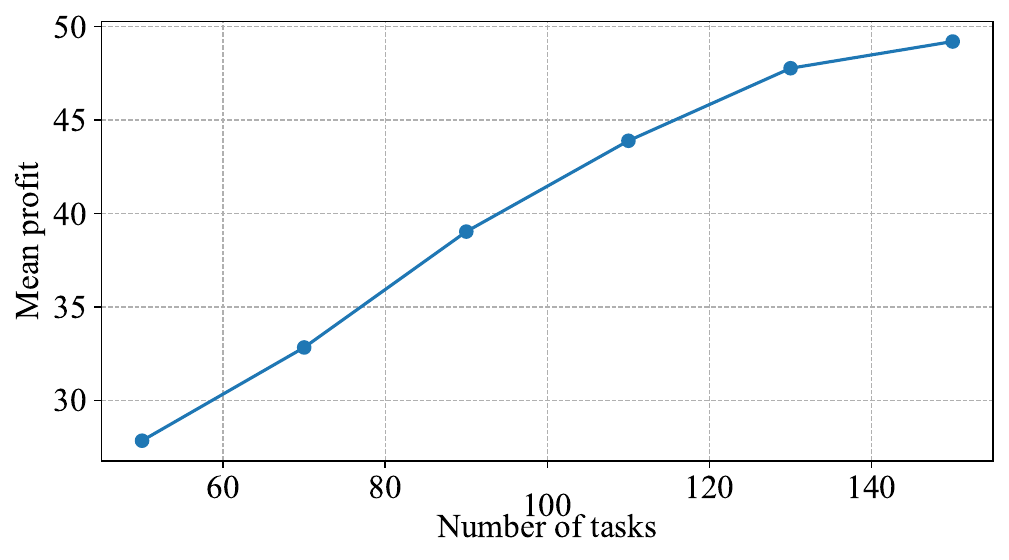}
	\caption{Mean profit $Z_{\text{S-ALNS}}$ as a function of the number of tasks with a fixed workforce size (10 workers) on map-based instances.}
	\label{fig:profit_vs_tasks_10v}
\end{figure}

To further isolate the effect of task volume, we additionally fix the workforce size at 10 and vary the number of tasks. Figure~\ref{fig:profit_vs_tasks_10v} shows that the mean profit increases steadily as more tasks become available, with diminishing marginal gains at larger $N$. This behavior is expected in prize-collecting routing: additional tasks provide more opportunities for profitable selections, while the finite time budgets and time windows eventually limit the achievable profit. Complete per-instance results for these configurations are reported in \ref{append3}.

\subsubsection{Effect of scenario-based lookahead}
\label{subsubsec:lookahead_effect}

To isolate the contribution of the scenario-based augmentation, we compare Scen-RH-ALNS with the myopic rolling-horizon variant introduced in Section~\ref{subsubsec:impl}. Myopic-RH uses the same HT-TOPTW snapshot model and ALNS components as Scen-RH-ALNS but sets the number of scenarios to one and does not include virtual tasks, so that each decision epoch is optimized solely based on currently revealed tasks.

On a collection of 31 map-based sets, Myopic-RH already produces solutions very close to the MIP incumbents obtained within the time limit. Adding scenario-based virtual tasks yields a further reduction in the average gap to the MIP incumbent of about $0.36$ percentage points. Scen-RH-ALNS improves the average gap on 22 of the 31 benchmark sets, and a paired $t$-test on the per-set differences indicates that this improvement, while numerically modest, is statistically significant ($t = 2.62$, two-sided $p \approx 0.009$, Cohen's $d \approx 0.47$). These results demonstrate that the scenario-based lookahead provides a systematic, albeit small, benefit over an already strong myopic rolling-horizon ALNS baseline.

Given that the residual gaps to the MIP incumbents on the map-based instances are only a few percentage points on average (Tables~\ref{tab:robust_dist} and~\ref{tab:scalability}), the incremental improvement due to scenarios is necessarily limited in absolute magnitude. Nevertheless, the consistent direction and statistical significance of the effect support the use of lightweight scenario-based augmentation in dynamic profitable routing policies.

A sensitivity analysis on the impact of parameters $S$ and $N_{\text{vir}}$ is presented in \ref{app:sensitivity_pilot}.

\subsubsection{Discussion}
\label{subsubsec:roadnet_discussion}

The map-based experiments indicate that Scen-RH-ALNS delivers stable and near-offline performance on heterogeneous DTOP-SC instances. Across all families in Table~\ref{tab:robust_dist}, the average gaps with respect to $Z_{MIP}$ remain within a moderate range, and in several configurations $\text{Gap}_{MIP}$ is close to zero or slightly negative, reflecting the ability of the online policy to match or occasionally outperform the MIP incumbent within the time limit. 

The standard deviations reported in Tables~\ref{tab:robust_dist} and~\ref{tab:scalability}, together with the per-instance gap distributions in Figures~\ref{fig:roadnet_gaps}, further indicate that the variability of Scen-RH-ALNS across instances is modest. Changes in service durations, time-window widths, profit ranges, and problem scale influence the absolute difficulty of the instances and the required time, but do not lead to instability or failure of the policy. Combined with the statistically significant improvement over Myopic-RH documented in Section~\ref{subsubsec:lookahead_effect}, these observations suggest that the proposed scenario-sampling rolling-horizon framework yields consistent solution quality for DTOP-SC on map-based instances while remaining computationally tractable.

\section{Conclusion}
\label{sec:conclusion}

This paper introduced the Dynamic Team Orienteering Problem in Spatial Crowdsourcing (DTOP-SC), addressing the challenge of coordinating heterogeneous fleets under dynamic demand and tight constraints. We formulated a static MIP model to clarify the problem structure and proposed the Scen-RH-ALNS policy. This framework combines rolling-horizon control with a scenario-sampling mechanism and aggregates worker-task assignment frequencies across scenarios to generate joint dispatch decisions in real time.

Computational experiments demonstrate the effectiveness and practicality of the proposed approach. On synthetic map-based DTOP-SC instances, Scen-RH-ALNS produces solutions whose profits are within a few percentage points of offline MIP incumbents across a wide range of instance families and scales, while average times remain moderate as problem size grows. On standard DTOP benchmark instances, Scen-RH-ALNS requires substantially smaller indicative runtimes than those reported for the multiple plan approaches on their respective tested platforms, while achieving average profits within 1--3\% of these methods, both in aggregate and when measured relative to CP-based offline solutions.

These results, together with the ablation against a myopic rolling-horizon variant, suggest that explicit scenario-based lookahead provides modest but consistent improvements and can help manage the trade-off between early commitment and flexibility. Future research avenues include integrating data-driven forecasting (e.g., machine learning) into the scenario generation layer, extending the model to handle stochastic travel times and service-level objectives, and developing adaptive mechanisms to tune algorithmic parameters online based on the system state. Another promising direction is to investigate deep reinforcement learning approaches that learn dispatch policies directly from data or simulation.

\section*{Acknowledgments}
This work was supported by the National Natural Science Foundation of China under Grants 72371175 and 71971148.

\section*{Declaration of generative AI and AI-assisted technologies in the manuscript preparation process}
During the preparation of this manuscript, the authors used ChatGPT in order to improve the readability and language of the text.
After using this tool, the authors reviewed and edited the content as needed and take full responsibility for the content of the published article.

\setcounter{table}{0}
\renewcommand{\thetable}{A\arabic{table}}

\clearpage

\appendix
\section{Additional computational results}

Table~\ref{tab:scalability_detailed} reports detailed scalability results for Scen-RH-ALNS on the map-based instances, including additional combinations of fleet size and task count. The performance trends are consistent with those observed in Table~\ref{tab:scalability} in the main text.

\begin{table}[htbp]
	\centering
	\caption{Scalability of Scen-RH-ALNS on additional map-based instances. S-ALNS denotes Scen-RH-ALNS.}
	\scriptsize
	\setlength{\tabcolsep}{4pt} 
	\renewcommand{\arraystretch}{1}
\begin{tabular}{m{1.2cm}<{\centering} m{1.8cm}<{\centering} m{2.0cm}<{\centering} m{1.8cm}<{\centering} m{1.5cm}<{\centering} m{3.0cm}<{\centering}}
		\toprule
		Dataset & Avg. Z$_{MIP}$& Avg. $ Z_{\text{S-ALNS}}$ & Gap$_{MIP}$& SD& Mean Time (s)\\
		\midrule
		5v/75t   & 25.71  & 23.79  & 7.48\%    & 2.30  & 4.90  \\
		5v/100t  & 30.77  & 29.56  & 3.93\%    & 4.21  & 13.47  \\
		7v/105t  & 36.06  & 34.22  & 5.11\%    & 2.59  & 12.88  \\
		7v/140t  & 40.69  & 40.25  & 1.08\%    & 5.24  & 19.02  \\
		9v/135t  & 43.15  & 42.39  & 1.77\%    & 5.83  & 22.40  \\
		9v/180t  & 37.57  & 47.18  & -25.57\%  & 7.33  & 19.67  \\
		10v/50t  & 28.33  & 27.84  & 1.74\%    & 2.67  & 10.25  \\
		10v/70t  & 33.84  & 32.83  & 3.00\%    & 3.91  & 12.47  \\
		10v/90t  & 40.82  & 39.03  & 4.38\%    & 2.70  & 17.66  \\
		10v/110t & 45.21  & 43.89  & 2.92\%    & 3.69  & 25.71  \\
		10v/130t & 47.11  & 47.77  & -1.39\%   & 6.00  & 27.26  \\
		10v/150t & 47.33  & 49.20  & -3.96\%   & 4.49  & 25.61  \\
		11v/165t & 40.52  & 56.66  & -39.83\%  & 4.96  & 38.87  \\
		11v/220t & 38.05  & 60.47  & -58.93\%  & 7.81  & 53.82  \\
		13v/195t & 40.77  & 65.80  & -61.39\%  & 5.35  & 75.48  \\
		13v/260t & 29.52  & 71.25  & -141.36\% & 5.64  & 94.44  \\
		15v/225t & 31.64  & 69.28  & -118.96\% & 6.74  & 61.96  \\
		15v/300t & 22.29  & 82.67  & -270.93\% & 5.28  & 112.04  \\
		\bottomrule
	\end{tabular}%
	\label{tab:scalability_detailed}%
\end{table}%

%
%
\clearpage

\section{Detail result on DTOP Benchmark}\label{append1}
\begin{table}[htbp]
	\centering
	\caption{Instance p1. Bold values indicate the new best-known online solution.}
	\fontsize{3.5pt}{5pt}\selectfont
	\resizebox{\linewidth}{!}{%
		\setlength{\tabcolsep}{2pt} 
		\renewcommand{\arraystretch}{1}
		\begin{tabular}{cccccccccccccccccccccccc}
			\toprule
			\multirow{2}[4]{*}{Problem} & \multirow{2}[4]{*}{ZCP} & \multicolumn{2}{c}{Scen-RH-ALNS} & \multicolumn{2}{c}{MPAd} & \multicolumn{2}{c}{MPAc} & \multirow{2}[4]{*}{Problem} & \multirow{2}[4]{*}{ZCP} & \multicolumn{2}{c}{Scen-RH-ALNS} & \multicolumn{2}{c}{MPAd} & \multicolumn{2}{c}{MPAc} & \multirow{2}[4]{*}{Problem} & \multirow{2}[4]{*}{ZCP} & \multicolumn{2}{c}{Scen-RH-ALNS} & \multicolumn{2}{c}{MPAd} & \multicolumn{2}{c}{MPAc} \\
			\cmidrule{3-8}\cmidrule{11-16}\cmidrule{19-24}          &       & Z     & Time   & Z     & Time   & Z     & Time   &       &       & Z     & Time   & Z     & Time   & Z     & Time   &       &       & Z     & Time   & Z     & Time   & Z     & Time \\
			\midrule
			dh.p1.2.a & 0     & 0     & 0.00  & 0     & 2.32  & 0     & 2.04  & dm.p1.2.a & 0     & 0     & 0.00  & 0     & 0.58  & 0     & 0.79  & dw.p1.2.a & 0     & 0     & 0.00  & 0     & 6.17  & 0     & 6.36  \\
			dh.p1.2.b & 0     & 0     & 0.00  & 0     & 3.60  & 0     & 0.75  & dm.p1.2.b & 15    & 0     & 0.00  & 0     & 1.03  & 0     & 1.38  & dw.p1.2.b & 5     & 5     & 1.82  & 5     & 9.34  & 5     & 6.22  \\
			dh.p1.2.c & 0     & 0     & 0.00  & 0     & 1.29  & 0     & 1.35  & dm.p1.2.c & 15    & 0     & 0.00  & 0     & 1.15  & 0     & 1.54  & dw.p1.2.c & 20    & 20    & 0.11  & 20    & 1.28  & 20    & 1.09  \\
			dh.p1.2.d & 5     & 0     & 0.00  & 0     & 1.78  & 0     & 1.68  & dm.p1.2.d & 10    & 5     & 1.83  & 5     & 1.56  & 5     & 2.09  & dw.p1.2.d & 30    & 25    & 0.06  & 25    & 1.55  & 25    & 1.57  \\
			dh.p1.2.e & 5     & 5     & 1.98  & 5     & 1.74  & 5     & 1.70  & dm.p1.2.e & 40    & 15    & 0.11  & 15    & 1.02  & 15    & 4.39  & dw.p1.2.e & 45    & \textbf{35} & 0.07  & 25    & 1.27  & 25    & 1.34  \\
			dh.p1.2.f & 25    & 15    & 0.06  & 15    & 0.65  & 15    & 0.67  & dm.p1.2.f & 45    & 25    & 0.06  & 25    & 0.98  & 25    & 5.55  & dw.p1.2.f & 70    & 60    & 0.09  & 60    & 1.32  & 60    & 1.44  \\
\bottomrule
	\end{tabular}}%
	\label{tab:p1}%
\end{table}%

\begin{table}[htbp]
	\centering
	\caption{Instance p2. Bold values indicate the new best-known online solution.}
	\fontsize{3.5pt}{5pt}\selectfont
	\resizebox{\linewidth}{!}{%
		\setlength{\tabcolsep}{2pt} 
		\renewcommand{\arraystretch}{1}
		\begin{tabular}{cccccccccccccccccccccccc}
			\toprule
			\multirow{2}[4]{*}{Problem} & \multirow{2}[4]{*}{ZCP} & \multicolumn{2}{c}{Scen-RH-ALNS} & \multicolumn{2}{c}{MPAd} & \multicolumn{2}{c}{MPAc} & \multirow{2}[4]{*}{Problem} & \multirow{2}[4]{*}{ZCP} & \multicolumn{2}{c}{Scen-RH-ALNS} & \multicolumn{2}{c}{MPAd} & \multicolumn{2}{c}{MPAc} & \multirow{2}[4]{*}{Problem} & \multirow{2}[4]{*}{ZCP} & \multicolumn{2}{c}{Scen-RH-ALNS} & \multicolumn{2}{c}{MPAd} & \multicolumn{2}{c}{MPAc} \\
			\cmidrule{3-8}\cmidrule{11-16}\cmidrule{19-24}          &       & Z     & Time   & Z     & Time   & Z     & Time   &       &       & Z     & Time   & Z     & Time   & Z     & Time   &       &       & Z     & Time   & Z     & Time   & Z     & Time \\
			\midrule
			dh.p2.2.a & 10    & 0     & 0.00  & 0     & 0.66  & 0     & 0.76  & dm.p2.2.a & 55    & 35    & 0.02  & 35    & 0.13  & 35    & 0.13  & dw.p2.2.a & 90    & 70    & 0.10  & 70    & 0.65  & 70    & 3.19  \\
			dh.p2.2.b & 25    & 0     & 0.00  & 0     & 0.70  & 0     & 0.76  & dm.p2.2.b & 90    & 35    & 0.05  & 35    & 0.12  & 35    & 0.12  & dw.p2.2.b & 120   & 95    & 0.14  & 95    & 0.54  & 95    & 2.90  \\
			dh.p2.2.c & 25    & 0     & 0.00  & 0     & 0.67  & 0     & 0.80  & dm.p2.2.c & 110   & 55    & 0.04  & 55    & 0.13  & 55    & 0.13  & dw.p2.2.c & 130   & \textbf{120} & 0.15  & 110   & 0.55  & 110   & 2.76  \\
			dh.p2.2.d & 45    & 20    & 0.05  & 20    & 0.58  & 20    & 0.59  & dm.p2.2.d & 115   & \textbf{75} & 0.10  & 65    & 0.14  & 65    & 0.15  & dw.p2.2.d & 160   & 125   & 0.21  & 135   & 0.34  & 130   & 1.70  \\
			dh.p2.2.e & 10    & 10    & 0.03  & 10    & 0.41  & 10    & 0.45  & dm.p2.2.e & 110   & 40    & 0.04  & 40    & 0.23  & 40    & 0.23  & dw.p2.2.e & 180   & 155   & 0.25  & 160   & 0.75  & 160   & 3.72  \\
\bottomrule
	\end{tabular}}%
	\label{tab:p2}%
\end{table}%

\begin{table}[htbp]
	\centering
	\caption{Instance p3. Bold values indicate the new best-known online solution.}
	\fontsize{3.5pt}{5pt}\selectfont
	\resizebox{\linewidth}{!}{%
		\setlength{\tabcolsep}{2pt} 
		\renewcommand{\arraystretch}{1}
		\begin{tabular}{cccccccccccccccccccccccc}
			\toprule
			\multirow{2}[4]{*}{Problem} & \multirow{2}[4]{*}{ZCP} & \multicolumn{2}{c}{Scen-RH-ALNS} & \multicolumn{2}{c}{MPAd} & \multicolumn{2}{c}{MPAc} & \multirow{2}[4]{*}{Problem} & \multirow{2}[4]{*}{ZCP} & \multicolumn{2}{c}{Scen-RH-ALNS} & \multicolumn{2}{c}{MPAd} & \multicolumn{2}{c}{MPAc} & \multirow{2}[4]{*}{Problem} & \multirow{2}[4]{*}{ZCP} & \multicolumn{2}{c}{Scen-RH-ALNS} & \multicolumn{2}{c}{MPAd} & \multicolumn{2}{c}{MPAc} \\
			\cmidrule{3-8}\cmidrule{11-16}\cmidrule{19-24}          &       & Z     & Time   & Z     & Time   & Z     & Time   &       &       & Z     & Time   & Z     & Time   & Z     & Time   &       &       & Z     & Time   & Z     & Time   & Z     & Time \\
			\midrule
			dh.p3.2.a & 10    & 10    & 0.02  & 10    & 4.05  & 10    & 3.78  & dm.p3.2.a & 70    & 10    & 0.02  & 10    & 3.13  & 10    & 1.60  & dw.p3.2.a & 40    & 40    & 0.06  & 40    & 1.77  & 40    & 2.02  \\
			dh.p3.2.b & 80    & 30    & 0.05  & 30    & 2.52  & 30    & 2.48  & dm.p3.2.b & 120   & 60    & 0.07  & 60    & 0.83  & 60    & 0.92  & dw.p3.2.b & 130   & 80    & 0.10  & 80    & 2.13  & 80    & 2.04  \\
			dh.p3.2.c & 50    & 40    & 0.07  & 40    & 5.30  & 40    & 3.94  & dm.p3.2.c & 160   & 70    & 0.10  & 70    & 1.38  & 70    & 1.36  & dw.p3.2.c & 160   & 130   & 0.11  & 130   & 1.50  & 120   & 1.54  \\
			dh.p3.2.d & 60    & 0     & 0.00  & 0     & 1.13  & 0     & 1.17  & dm.p3.2.d & 150   & 80    & 0.15  & 80    & 0.87  & 80    & 1.13  & dw.p3.2.d & 200   & 170   & 0.11  & 160   & 1.49  & 170   & 4.92  \\
			dh.p3.2.e & 70    & 20    & 0.03  & 20    & 0.97  & 20    & 0.97  & dm.p3.2.e & 190   & 120   & 0.19  & 120   & 1.12  & 120   & 1.38  & dw.p3.2.e & 250   & 200   & 0.20  & 220   & 4.76  & 220   & 10.69  \\
			dh.p3.2.f & 100   & 30    & 0.06  & 30    & 1.51  & 30    & 1.52  & dm.p3.2.f & 230   & 80    & 0.15  & 80    & 1.21  & 120   & 1.36  & dw.p3.2.f & 290   & 200   & 0.30  & 230   & 9.25  & 230   & 3.14  \\
\bottomrule
	\end{tabular}}%
	\label{tab:p3}%
\end{table}%

\begin{table}[htbp]
	\centering
	\caption{Instance p4. Bold values indicate the new best-known online solution.}
	\fontsize{3.5pt}{5pt}\selectfont
	\resizebox{\linewidth}{!}{%
		\setlength{\tabcolsep}{2pt} 
		\renewcommand{\arraystretch}{1}
		\begin{tabular}{cccccccccccccccccccccccc}
			\toprule
			\multirow{2}[4]{*}{Problem} & \multirow{2}[4]{*}{ZCP} & \multicolumn{2}{c}{Scen-RH-ALNS} & \multicolumn{2}{c}{MPAd} & \multicolumn{2}{c}{MPAc} & \multirow{2}[4]{*}{Problem} & \multirow{2}[4]{*}{ZCP} & \multicolumn{2}{c}{Scen-RH-ALNS} & \multicolumn{2}{c}{MPAd} & \multicolumn{2}{c}{MPAc} & \multirow{2}[4]{*}{Problem} & \multirow{2}[4]{*}{ZCP} & \multicolumn{2}{c}{Scen-RH-ALNS} & \multicolumn{2}{c}{MPAd} & \multicolumn{2}{c}{MPAc} \\
			\cmidrule{3-8}\cmidrule{11-16}\cmidrule{19-24}          &       & Z     & Time   & Z     & Time   & Z     & Time   &       &       & Z     & Time   & Z     & Time   & Z     & Time   &       &       & Z     & Time   & Z     & Time   & Z     & Time \\
			\midrule
			dh.p4.2.a & 84    & 28    & 0.02  & 28    & 48.96  & 28    & 48.44  & dm.p4.2.a & 146   & 70    & 0.07  & 70    & 155.59  & 70    & 157.16  & dw.p4.2.a & 143   & 102   & 0.16  & 125   & 215.65  & 125   & 140.38  \\
			dh.p4.2.b & 127   & 86    & 0.07  & 86    & 47.90  & 86    & 47.76  & dm.p4.2.b & 215   & 111   & 0.11  & 111   & 179.12  & 111   & 174.52  & dw.p4.2.b & 281   & 202   & 0.65  & 234   & 326.71  & 230   & 303.07  \\
			dh.p4.2.c & 156   & 101   & 0.12  & 101   & 207.78  & 101   & 208.02  & dm.p4.2.c & 273   & 175   & 0.28  & 171   & 294.53  & 201   & 294.04  & dw.p4.2.c & 380   & 331   & 1.87  & 351   & 426.57  & 351   & 411.73  \\
			dh.p4.2.d & 180   & 120   & 0.11  & 120   & 289.81  & 120   & 288.81  & dm.p4.2.d & 335   & \textbf{286} & 0.57  & 272   & 370.46  & 272   & 370.37  & dw.p4.2.d & 471   & 269   & 1.67  & 433   & 554.46  & 424   & 498.82  \\
			dh.p4.2.e & 206   & 69    & 0.14  & 87    & 392.12  & 87    & 387.05  & dm.p4.2.e & 399   & \textbf{338} & 1.67  & 306   & 398.05  & 306   & 377.84  & dw.p4.2.e & 521   & 362   & 2.54  & 460   & 591.32  & 468   & 581.31  \\
\bottomrule
	\end{tabular}}%
	\label{tab:p4}%
\end{table}%

\begin{table}[htbp]
	\centering
	\caption{Instance p5. Bold values indicate the new best-known online solution.}
	\fontsize{3.5pt}{5pt}\selectfont
	\resizebox{\linewidth}{!}{%
		\setlength{\tabcolsep}{2pt} 
		\renewcommand{\arraystretch}{1}
		\begin{tabular}{cccccccccccccccccccccccc}
			\toprule
			\multirow{2}[4]{*}{Problem} & \multirow{2}[4]{*}{ZCP} & \multicolumn{2}{c}{Scen-RH-ALNS} & \multicolumn{2}{c}{MPAd} & \multicolumn{2}{c}{MPAc} & \multirow{2}[4]{*}{Problem} & \multirow{2}[4]{*}{ZCP} & \multicolumn{2}{c}{Scen-RH-ALNS} & \multicolumn{2}{c}{MPAd} & \multicolumn{2}{c}{MPAc} & \multirow{2}[4]{*}{Problem} & \multirow{2}[4]{*}{ZCP} & \multicolumn{2}{c}{Scen-RH-ALNS} & \multicolumn{2}{c}{MPAd} & \multicolumn{2}{c}{MPAc} \\
			\cmidrule{3-8}\cmidrule{11-16}\cmidrule{19-24}          &       & Z     & Time   & Z     & Time   & Z     & Time   &       &       & Z     & Time   & Z     & Time   & Z     & Time   &       &       & Z     & Time   & Z     & Time   & Z     & Time \\
			\midrule
			    dh.p5.2.a & 0     & 0     & 0.00  & 0     & 31.96  & 0     & 31.13  & dm.p5.2.a & 0     & 0     & 0.00  & 0     & 65.00  & 0     & 62.27  & dw.p5.2.a & 0     & 0     & 0.00  & 0     & 81.45  & 0     & 74.94  \\
			dh.p5.2.b & 10    & 5     & 0.02  & 5     & 54.64  & 5     & 53.95  & dm.p5.2.b & 10    & 0     & 0.00  & 0     & 109.97  & 0     & 112.73  & dw.p5.2.b & 20    & 15    & 0.04  & 15    & 49.39  & 15    & 49.05  \\
			dh.p5.2.c & 10    & 0     & 0.00  & 0     & 52.27  & 0     & 52.07  & dm.p5.2.c & 40    & 30    & 0.03  & 30    & 61.40  & 30    & 58.91  & dw.p5.2.c & 45    & 40    & 0.05  & 40    & 51.61  & 40    & 48.93  \\
			dh.p5.2.d & 40    & 5     & 0.02  & 5     & 121.83  & 5     & 121.63  & dm.p5.2.d & 70    & 45    & 0.04  & 45    & 55.58  & 45    & 71.49  & dw.p5.2.d & 80    & 70    & 0.09  & 70    & 28.39  & 70    & 27.76  \\
			dh.p5.2.e & 45    & 10    & 0.03  & 10    & 119.41  & 10    & 119.36  & dm.p5.2.e & 130   & 80    & 0.08  & 80    & 49.28  & 80    & 49.87  & dw.p5.2.e & 170   & 120   & 0.16  & 140   & 37.61  & 140   & 30.32  \\
\bottomrule
	\end{tabular}}%
	\label{tab:p5}%
\end{table}%

\begin{table}[htbp]
	\centering
	\caption{Instance p6. Bold values indicate the new best-known online solution.}
	\fontsize{3.5pt}{5pt}\selectfont
	\resizebox{\linewidth}{!}{%
		\setlength{\tabcolsep}{2pt} 
		\renewcommand{\arraystretch}{1}
		\begin{tabular}{cccccccccccccccccccccccc}
			\toprule
			\multirow{2}[4]{*}{Problem} & \multirow{2}[4]{*}{ZCP} & \multicolumn{2}{c}{Scen-RH-ALNS} & \multicolumn{2}{c}{MPAd} & \multicolumn{2}{c}{MPAc} & \multirow{2}[4]{*}{Problem} & \multirow{2}[4]{*}{ZCP} & \multicolumn{2}{c}{Scen-RH-ALNS} & \multicolumn{2}{c}{MPAd} & \multicolumn{2}{c}{MPAc} & \multirow{2}[4]{*}{Problem} & \multirow{2}[4]{*}{ZCP} & \multicolumn{2}{c}{Scen-RH-ALNS} & \multicolumn{2}{c}{MPAd} & \multicolumn{2}{c}{MPAc} \\
			\cmidrule{3-8}\cmidrule{11-16}\cmidrule{19-24}          &       & Z     & Time   & Z     & Time   & Z     & Time   &       &       & Z     & Time   & Z     & Time   & Z     & Time   &       &       & Z     & Time   & Z     & Time   & Z     & Time \\
			\midrule
			    dh.p6.2.a & 0     & 0     & 0.00  & 0     & 5.96  & 0     & 4.26  & dm.p6.2.a & 0     & 0     & 0.00  & 0     & 3.62  & 0     & 3.25  & dw.p6.2.a & 0     & 0     & 0.00  & 0     & 2.58  & 0     & 2.61  \\
			dh.p6.2.b & 0     & 0     & 0.00  & 0     & 4.09  & 0     & 15.57  & dm.p6.2.b & 0     & 0     & 0.00  & 0     & 3.11  & 0     & 3.37  & dw.p6.2.b & 0     & 0     & 0.00  & 0     & 2.52  & 0     & 2.70  \\
			dh.p6.2.c & 0     & 0     & 0.00  & 0     & 10.20  & 0     & 5.23  & dm.p6.2.c & 0     & 0     & 0.00  & 0     & 3.25  & 0     & 13.74  & dw.p6.2.c & 0     & 0     & 0.00  & 0     & 2.67  & 0     & 2.80  \\
			dh.p6.2.d & 108   & 54    & 0.06  & 54    & 10.62  & 54    & 4.13  & dm.p6.2.d & 174   & \textbf{162} & 0.13  & 156   & 30.59  & 156   & 32.10  & dw.p6.2.d & 186   & 156   & 0.24  & 156   & 11.54  & 156   & 11.50  \\
			dh.p6.2.e & 108   & 72    & 0.09  & 72    & 18.27  & 72    & 22.40  & dm.p6.2.e & 288   & 210   & 0.36  & 216   & 50.02  & 216   & 37.15  & dw.p6.2.e & 348   & 306   & 0.75  & 318   & 23.51  & 318   & 23.19  \\
\bottomrule
	\end{tabular}}%
	\label{tab:p6}%
\end{table}%

\begin{table}[htbp]
	\centering
	\caption{Instance p7. Bold values indicate the new best-known online solution.}
	\fontsize{3.5pt}{5pt}\selectfont
	\resizebox{\linewidth}{!}{%
		\setlength{\tabcolsep}{2pt} 
		\renewcommand{\arraystretch}{1}
		\begin{tabular}{cccccccccccccccccccccccc}
			\toprule
			\multirow{2}[4]{*}{Problem} & \multirow{2}[4]{*}{ZCP} & \multicolumn{2}{c}{Scen-RH-ALNS} & \multicolumn{2}{c}{MPAd} & \multicolumn{2}{c}{MPAc} & \multirow{2}[4]{*}{Problem} & \multirow{2}[4]{*}{ZCP} & \multicolumn{2}{c}{Scen-RH-ALNS} & \multicolumn{2}{c}{MPAd} & \multicolumn{2}{c}{MPAc} & \multirow{2}[4]{*}{Problem} & \multirow{2}[4]{*}{ZCP} & \multicolumn{2}{c}{Scen-RH-ALNS} & \multicolumn{2}{c}{MPAd} & \multicolumn{2}{c}{MPAc} \\
			\cmidrule{3-8}\cmidrule{11-16}\cmidrule{19-24}          &       & Z     & Time   & Z     & Time   & Z     & Time   &       &       & Z     & Time   & Z     & Time   & Z     & Time   &       &       & Z     & Time   & Z     & Time   & Z     & Time \\
			\midrule
			    dh.p7.2.a & 0     & 0     & 0.00  & 0     & 689.22  & 0     & 568.51  & dm.p7.2.a & 0     & 0     & 0.01  & 0     & 1191.62  & 0     & 1192.95  & dw.p7.2.a & 30    & 14    & 0.01  & 14    & 861.35  & 14    & 842.98  \\
			dh.p7.2.b & 16    & 16    & 0.02  & 16    & 958.83  & 16    & 1020.21  & dm.p7.2.b & 30    & 14    & 0.02  & 14    & 1700.27  & 14    & 1700.89  & dw.p7.2.b & 64    & 29    & 0.03  & 29    & 890.26  & 29    & 879.00  \\
			dh.p7.2.c & 23    & 0     & 0.01  & 0     & 1154.22  & 0     & 1170.17  & dm.p7.2.c & 86    & 58    & 0.09  & 58    & 1050.35  & 58    & 1053.48  & dw.p7.2.c & 101   & 48    & 0.08  & 48    & 413.73  & 48    & 412.51  \\
			dh.p7.2.d & 49    & 0     & 0.01  & 0     & 1290.80  & 0     & 1227.60  & dm.p7.2.d & 147   & 104   & 0.16  & 104   & 572.46  & 104   & 572.63  & dw.p7.2.d & 190   & 158   & 0.18  & 179   & 235.46  & 179   & 240.13  \\
			dh.p7.2.e & 90    & 55    & 0.03  & 60    & 106.29  & 60    & 89.56  & dm.p7.2.e & 213   & 158   & 0.24  & 197   & 506.66  & 197   & 516.67  & dw.p7.2.e & 267   & 226   & 0.39  & 250   & 300.39  & 250   & 294.50  \\
\bottomrule
	\end{tabular}}%
	\label{tab:p7}%
\end{table}%

\clearpage

\subsection{Pilot sensitivity to the numbers of scenarios and virtual tasks}
\label{app:sensitivity_pilot}

We conducted a small pilot sweep on a representative Base subset (10 instances) to assess the sensitivity of Scen-RH-ALNS to the number of sampled scenarios $S$ and the number of virtual tasks per scenario $N_{\text{vir}}$. Across the tested settings, profits showed substantial overlap and no clear monotone trend in either $S$ or $N_{\text{vir}}$, while decision time increased consistently as $S$ or $N_{\text{vir}}$ grew, suggesting diminishing returns in solution quality relative to computational effort. These observations motivate focusing the main ablation on the two semantically distinct configurations Myopic-RH ($S=1,\,N_{\text{vir}}=0$) and Scen-RH-ALNS (default $S=15,\,N_{\text{vir}}=5$), with statistical conclusions drawn from the full benchmark collection.

\begin{table}[!htbp]
	\centering
	\caption{Pilot sensitivity of Scen-RH-ALNS to $S$ on a Base subset (10 instances): mean profit $Avg. Z_{\text{S-ALNS}}$, SD, and mean time.}	
	\label{tab:sensitivity_pilot_S}
	\scriptsize
	\setlength{\tabcolsep}{3pt}
	\renewcommand{\arraystretch}{1}
	\begin{tabular}{cccc}
		\toprule Setting & $Avg. Z_{\text{S-ALNS}}$ & SD & Mean Time (s)\\
		\midrule 
		$S=1$ & 40.534 & 4.40 & 9.81\\ 
		$S=3$ & 40.416 & 4.65 & 10.73\\ 
		$S=5$ & 40.788 & 4.23 & 11.21\\ 
		$S=7$ & 40.374 & 4.66 & 12.85\\ 
		$S=9$ & 40.590 & 4.65 & 14.78\\ 
		$S=11$ & 40.740 & 4.63 & 16.05\\ 
		$S=13$ & 40.474 & 4.62 & 17.79\\ 
		$S=15$ & 40.562 & 4.52 & 19.13\\
		
		\bottomrule
	\end{tabular}
\end{table}

\begin{table}[!htbp]
	\centering
	\caption{Pilot sensitivity of Scen-RH-ALNS to $N_{\text{vir}}$ on a Base subset (10 instances): mean profit $Avg. Z_{\text{S-ALNS}}$, SD, and mean time.}	
	\label{tab:sensitivity_pilot_N}
	\scriptsize
	\setlength{\tabcolsep}{3pt}
	\renewcommand{\arraystretch}{1}
	\begin{tabular}{cccc}
		\toprule 
		Setting & $Avg. Z_{\text{S-ALNS}}$ & SD & Mean Time (s)\\
		\midrule
		$N_{\text{vir}}=1$ & 40.340 & 4.38 & 15.28\\
		$N_{\text{vir}}=2$ & 40.462 & 4.35 & 16.24\\ 
		$N_{\text{vir}}=3$ & 40.652 & 4.41 & 16.84\\ 
		$N_{\text{vir}}=4$ & 40.546 & 4.70 & 18.40\\ 
		$N_{\text{vir}}=5$ & 40.626 & 4.37 & 19.49\\ 
		$N_{\text{vir}}=6$ & 40.566 & 4.73 & 20.95\\ 
		$N_{\text{vir}}=7$ & 40.358 & 4.36 & 22.11\\ 
		$N_{\text{vir}}=8$ & 40.546 & 4.48 & 23.50\\ 
		$N_{\text{vir}}=9$ & 40.702 & 4.61 & 24.84\\ 
		$N_{\text{vir}}=10$ & 40.692 & 4.38 & 26.67\\
		\bottomrule
	\end{tabular}
\end{table}

\clearpage

\section{Detail result on road network instances}\label{append2}

In addition to the instance-wise gap $\text{Gap}_{\text{MIP}} = 100\% \times (Z_{\text{MIP}}-Z_{\text{S-ALNS}})/Z_{\text{MIP}}$, we report an aggregate gap that summarizes performance:
\[
\text{AggGap}_{\text{MIP}} = 100\% \times
\frac{\sum Z_{\text{MIP}} - \sum Z_{\text{S-ALNS}}}
{\sum Z_{\text{MIP}}}.
\]

\begin{table}[!htbp]
	\centering
	\caption{Instance rn.short.}
	\label{tab:short}
	\scriptsize
	\renewcommand{\arraystretch}{0.8}
	\begin{threeparttable}
\begin{tabular}{ccccc}
			\toprule 
			\makecell[c]{Instance\\(idx)} & $Z_{\text{S-ALNS}}$ & $Z_{MIP}$ & Time (s) & \makecell[c]{Gap to\\$Z_{MIP}$} \\
			\midrule
			rn.short.0 & 54.16  & 53.42  & 73.45  & -1.39\% \\
			rn.short.1 & 53.30  & 53.14  & 49.54  & -0.30\% \\
			rn.short.2 & 60.24  & 59.88  & 95.22  & -0.60\% \\
			rn.short.3 & 43.56  & 42.38  & 37.94  & -2.78\% \\
			rn.short.4 & 58.30  & 57.44  & 172.33 & -1.50\% \\
			rn.short.5 & 49.22  & 47.26  & 69.06  & -4.15\% \\
			rn.short.6 & 58.54  & 59.26  & 173.81 & 1.21\%  \\
			rn.short.7 & 57.14  & 55.62  & 106.68 & -2.73\% \\
			rn.short.8 & 48.60  & 47.86  & 94.71  & -1.55\% \\
			rn.short.9 & 59.92  & 59.88  & 142.07 & -0.07\% \\
			\midrule
			Mean  & 54.30  & 53.61  & 101.48  & -1.38\% \\
		\end{tabular}%
		
		\begin{tablenotes}
			\scriptsize
			\item \textit{Note:} Aggregate gap to $Z_{MIP}$ is -1.28\%.
		\end{tablenotes}
	\end{threeparttable}
\end{table}%

\begin{table}[!htbp]
	\centering
	\caption{Instance rn.long.}
	\label{tab:long}
	\scriptsize
	\renewcommand{\arraystretch}{0.8}
	\begin{threeparttable}
		\begin{tabular}{ccccc}
			\toprule 
			\makecell[c]{Instance\\(idx)} & $Z_{\text{S-ALNS}}$ & $Z_{MIP}$ & Time (s) & \makecell[c]{Gap to\\$Z_{MIP}$} \\
			\midrule
			rn.long.0 & 22.68  & 26.08  & 3.40  & 13.04\% \\
			rn.long.1 & 23.70  & 25.68  & 3.33  & 7.71\%  \\
			rn.long.2 & 29.80  & 30.20  & 4.46  & 1.32\%  \\
			rn.long.3 & 32.26  & 34.86  & 6.93  & 7.46\%  \\
			rn.long.4 & 28.48  & 28.20  & 3.81  & -0.99\% \\
			rn.long.5 & 26.58  & 28.36  & 3.91  & 6.28\%  \\
			rn.long.6 & 30.44  & 32.00  & 6.89  & 4.88\%  \\
			rn.long.7 & 28.00  & 27.70  & 7.25  & -1.08\% \\
			rn.long.8 & 22.74  & 23.98  & 3.05  & 5.17\%  \\
			rn.long.9 & 27.56  & 27.20  & 3.42  & -1.32\% \\
			\midrule
			Mean  & 27.22  & 28.43  & 4.65  & 4.25\% \\
		\end{tabular}%
		
		\begin{tablenotes}
			\scriptsize
			\item \textit{Note:} Aggregate gap to $Z_{MIP}$ is 4.23\%.
		\end{tablenotes}
	\end{threeparttable}
\end{table}
\begin{table}[!htbp]
\centering
\caption{Instance rn.tight.}
\label{tab:tight}
\scriptsize
\renewcommand{\arraystretch}{0.8}
\begin{threeparttable}
	\begin{tabular}{ccccc}
		\toprule 
		\makecell[c]{Instance\\(idx)} & $Z_{\text{S-ALNS}}$ & $Z_{MIP}$ & Time (s) & \makecell[c]{Gap to\\$Z_{MIP}$} \\
		\midrule
		rn.tight.0 & 42.80  & 45.10  & 14.78  & 5.10\% \\
		rn.tight.1 & 38.46  & 39.94  & 15.82  & 3.71\% \\
		rn.tight.2 & 34.14  & 35.88  & 8.28  & 4.85\% \\
		rn.tight.3 & 42.18  & 43.36  & 20.88  & 2.72\% \\
		rn.tight.4 & 46.64  & 47.64  & 15.09  & 2.10\% \\
		rn.tight.5 & 42.70  & 41.70  & 25.99  & -2.40\% \\
		rn.tight.6 & 43.98  & 46.34  & 15.04  & 5.09\% \\
		rn.tight.7 & 42.26  & 43.66  & 9.70  & 3.21\% \\
		rn.tight.8 & 36.64  & 37.22  & 26.85  & 1.56\% \\
		rn.tight.9 & 46.00  & 46.56  & 24.01  & 1.20\% \\
		\midrule
		Mean  & 41.58  & 42.74  & 17.64  & 2.71\% \\
	\end{tabular}%
	
	\begin{tablenotes}
		\scriptsize
		\item \textit{Note:} Aggregate gap to $Z_{MIP}$ is 2.71\%.
	\end{tablenotes}
\end{threeparttable}
\end{table}
\begin{table}[!htbp]
\centering
\caption{Instance rn.loose.}
\label{tab:loose}
\scriptsize
\renewcommand{\arraystretch}{0.8}
\begin{threeparttable}
\begin{tabular}{ccccc}
	\toprule 
	\makecell[c]{Instance\\(idx)} & $Z_{\text{S-ALNS}}$ & $Z_{MIP}$ & Time (s) & \makecell[c]{Gap to\\$Z_{MIP}$} \\
	\midrule
	rn.loose.0 & 35.10  & 37.06  & 13.75  & 5.29\% \\
	rn.loose.1 & 44.88  & 42.80  & 38.41  & -4.86\% \\
	rn.loose.2 & 40.00  & 42.94  & 20.40  & 6.85\% \\
	rn.loose.3 & 41.12  & 42.94  & 32.51  & 4.24\% \\
	rn.loose.4 & 35.58  & 36.52  & 21.20  & 2.57\% \\
	rn.loose.5 & 39.14  & 39.74  & 21.63  & 1.51\% \\
	rn.loose.6 & 35.44  & 32.90  & 15.97  & -7.72\% \\
	rn.loose.7 & 38.90  & 39.38  & 32.01  & 1.22\% \\
	rn.loose.8 & 38.54  & 35.04  & 25.52  & -9.99\% \\
	rn.loose.9 & 34.14  & 36.50  & 19.40  & 6.47\% \\
	\midrule
	Mean  & 38.28  & 38.58  & 24.08  & 0.56\% \\
\end{tabular}%

\begin{tablenotes}
	\scriptsize
	\item \textit{Note:} Aggregate gap to $Z_{MIP}$ is 0.77\%.
\end{tablenotes}
\end{threeparttable}
\end{table}
\begin{table}[!htbp]
\centering
\caption{Instance rn.narrow.}
\label{tab:narrow}
\scriptsize
\renewcommand{\arraystretch}{0.8}
\begin{threeparttable}
\begin{tabular}{ccccc}
\toprule 
\makecell[c]{Instance\\(idx)} & $Z_{\text{S-ALNS}}$ & $Z_{MIP}$  & Time (s) & \makecell[c]{Gap to\\$Z_{MIP}$}  \\
\midrule
rn.narrow.0 & 53.45  & 53.70  & 32.78  & 0.47\%  \\
rn.narrow.1 & 43.55  & 45.40  & 26.14  & 4.07\% \\
rn.narrow.2 & 43.85  & 44.55  & 22.67  & 1.57\% \\
rn.narrow.3 & 53.15  & 54.20  & 21.96  & 1.94\% \\
rn.narrow.4 & 50.05  & 53.20  & 30.96  & 5.92\% \\
rn.narrow.5 & 48.80  & 50.50  & 11.89  & 3.37\% \\
rn.narrow.6 & 55.80  & 53.45  & 33.40  & -4.40\% \\
rn.narrow.7 & 45.00  & 44.40  & 11.41  & -1.35\% \\
rn.narrow.8 & 48.10  & 48.65  & 13.84  & 1.13\% \\
rn.narrow.9 & 39.00  & 40.65  & 9.11  & 4.06\% \\
\midrule
Mean  & 48.08  & 48.87  & 21.42  & 1.68\% \\
\end{tabular}%

\begin{tablenotes}
\scriptsize
\item \textit{Note:} Aggregate gap to $Z_{MIP}$ is 1.63\%.
\end{tablenotes}
\end{threeparttable}
\end{table}
\begin{table}[!htbp]
\centering
\caption{Instance rn.wide.}
\label{tab:wide}
\scriptsize
\renewcommand{\arraystretch}{0.8}
\begin{threeparttable}
\begin{tabular}{ccccc}
\toprule 
\makecell[c]{Instance\\(idx)} & $Z_{\text{S-ALNS}}$ & $Z_{MIP}$ & Time (s) & \makecell[c]{Gap to\\$Z_{MIP}$} \\
\midrule
rn.wide.0 & 37.64  & 38.91  & 13.73  & 3.26\% \\
rn.wide.1 & 40.04  & 42.20  & 14.72  & 5.12\% \\
rn.wide.2 & 30.75  & 31.64  & 7.42  & 2.81\% \\
rn.wide.3 & 40.23  & 43.49  & 16.69  & 7.50\% \\
rn.wide.4 & 31.25  & 33.27  & 15.02  & 6.07\% \\
rn.wide.5 & 35.18  & 39.06  & 11.85  & 9.93\% \\
rn.wide.6 & 34.68  & 35.83  & 11.47  & 3.21\% \\
rn.wide.7 & 32.33  & 34.41  & 14.66  & 6.04\% \\
rn.wide.8 & 36.77  & 37.94  & 8.96  & 3.08\% \\
rn.wide.9 & 37.89  & 38.15  & 28.49  & 0.68\% \\
\midrule
Mean  & 35.68  & 37.49  & 14.30  & 4.77\% \\
\end{tabular}%

\begin{tablenotes}
\scriptsize
\item \textit{Note:} Aggregate gap to $Z_{MIP}$ is 4.84\%.
\end{tablenotes}
\end{threeparttable}
\end{table}
\begin{table}[!htbp]
\centering
\caption{Instance rn.5v/50t.}
\label{tab:5v/50t}
\scriptsize
\renewcommand{\arraystretch}{0.8}
\begin{threeparttable}
\begin{tabular}{ccccc}
\toprule 
\makecell[c]{Instance\\(idx)} & $Z_{\text{S-ALNS}}$ & $Z_{MIP}$ & Time (s) & \makecell[c]{Gap to\\$Z_{MIP}$} \\
\midrule
rn.5v/50t.0 & 21.14  & 21.98  & 3.01  & 3.82\% \\
rn.5v/50t.1 & 22.84  & 23.56  & 3.01  & 3.06\% \\
rn.5v/50t.2 & 24.18  & 24.88  & 5.88  & 2.81\% \\
rn.5v/50t.3 & 21.24  & 23.26  & 1.54  & 8.68\% \\
rn.5v/50t.4 & 20.12  & 21.00  & 3.06  & 4.19\% \\
rn.5v/50t.5 & 23.26  & 25.46  & 9.26  & 8.64\% \\
rn.5v/50t.6 & 22.02  & 22.70  & 4.52  & 3.00\% \\
rn.5v/50t.7 & 19.02  & 21.30  & 1.81  & 10.70\% \\
rn.5v/50t.8 & 22.32  & 22.60  & 2.42  & 1.24\% \\
rn.5v/50t.9 & 25.24  & 26.62  & 5.94  & 5.18\% \\
\midrule
Mean  & 22.14  & 23.34  & 4.05  & 5.13\% \\
\end{tabular}%

\begin{tablenotes}
\scriptsize
\item \textit{Note:} Aggregate gap to $Z_{MIP}$ is 5.13\%.
\end{tablenotes}
\end{threeparttable}
\end{table}
\begin{table}[!htbp]
\centering
\caption{Instance rn.7v/70t.}
\label{tab:7v/70t}
\scriptsize
\renewcommand{\arraystretch}{0.8}
\begin{threeparttable}
\begin{tabular}{ccccc}
\toprule 
\makecell[c]{Instance\\(idx)} & $Z_{\text{S-ALNS}}$ & $Z_{MIP}$ & Time (s) & \makecell[c]{Gap to\\$Z_{MIP}$} \\
\midrule
rn.7v/70t.0 & 27.54  & 30.26  & 6.98  & 8.99\% \\
rn.7v/70t.1 & 27.38  & 29.20  & 4.76  & 6.23\% \\
rn.7v/70t.2 & 27.74  & 29.64  & 11.65  & 6.41\% \\
rn.7v/70t.3 & 25.58  & 26.96  & 4.82  & 5.12\% \\
rn.7v/70t.4 & 29.04  & 29.86  & 7.25  & 2.75\% \\
rn.7v/70t.5 & 23.94  & 25.84  & 5.64  & 7.35\% \\
rn.7v/70t.6 & 27.74  & 30.36  & 8.27  & 8.63\% \\
rn.7v/70t.7 & 31.00  & 32.88  & 16.50  & 5.72\% \\
rn.7v/70t.8 & 34.40  & 35.02  & 9.81  & 1.77\% \\
rn.7v/70t.9 & 30.18  & 32.58  & 14.04  & 7.37\% \\
\midrule
Mean  & 28.45  & 30.26  & 8.97  & 6.03\% \\
\end{tabular}%

\begin{tablenotes}
\scriptsize
\item \textit{Note:} Aggregate gap to $Z_{MIP}$ is 5.97\%.
\end{tablenotes}
\end{threeparttable}
\end{table}
\begin{table}[!htbp]
\centering
\caption{Instance rn.9v/90t.}
\label{tab:9v/90t}
\scriptsize
\renewcommand{\arraystretch}{0.8}
\begin{threeparttable}
\begin{tabular}{ccccc}
\toprule 
\makecell[c]{Instance\\(idx)} & $Z_{\text{S-ALNS}}$ & $Z_{MIP}$ & Time (s) & \makecell[c]{Gap to\\$Z_{MIP}$} \\
\midrule
rn.9v/90t.0 & 38.10  & 40.36  & 10.46  & 5.60\% \\
rn.9v/90t.1 & 43.26  & 45.40  & 26.05  & 4.71\% \\
rn.9v/90t.2 & 35.84  & 42.18  & 15.22  & 15.03\% \\
rn.9v/90t.3 & 39.84  & 42.66  & 12.27  & 6.61\% \\
rn.9v/90t.4 & 38.06  & 40.10  & 19.08  & 5.09\% \\
rn.9v/90t.5 & 32.70  & 33.84  & 5.36  & 3.37\% \\
rn.9v/90t.6 & 35.58  & 38.20  & 10.26  & 6.86\% \\
rn.9v/90t.7 & 40.30  & 40.70  & 23.61  & 0.98\% \\
rn.9v/90t.8 & 31.86  & 32.80  & 6.84  & 2.87\% \\
rn.9v/90t.9 & 31.00  & 33.48  & 10.08  & 7.41\% \\
\midrule
Mean  & 36.65  & 38.97  & 13.92  & 5.85\% \\
\end{tabular}%

\begin{tablenotes}
\scriptsize
\item \textit{Note:} Aggregate gap to $Z_{MIP}$ is 5.95\%.
\end{tablenotes}
\end{threeparttable}
\end{table}
\begin{table}[!htbp]
\centering
\caption{Instance rn.10v/100t. (Base)}
\label{tab:10v/100t}
\scriptsize
\renewcommand{\arraystretch}{0.8}
\begin{threeparttable}
\begin{tabular}{ccccc}
\toprule 
\makecell[c]{Instance\\(idx)} & $Z_{\text{S-ALNS}}$ & $Z_{MIP}$ & Time (s) & \makecell[c]{Gap to\\$Z_{MIP}$} \\
\midrule
rn.10v/100t.0 & 42.40  & 43.84  & 23.65  & 3.28\% \\
rn.10v/100t.1 & 40.08  & 44.10  & 10.80  & 9.12\% \\
rn.10v/100t.2 & 38.52  & 40.16  & 25.77  & 4.08\% \\
rn.10v/100t.3 & 38.62  & 42.66  & 17.53  & 9.47\% \\
rn.10v/100t.4 & 36.22  & 39.06  & 12.70  & 7.27\% \\
rn.10v/100t.5 & 47.06  & 49.56  & 26.28  & 5.04\% \\
rn.10v/100t.6 & 50.18  & 51.56  & 23.79  & 2.68\% \\
rn.10v/100t.7 & 38.26  & 39.24  & 8.61  & 2.50\% \\
rn.10v/100t.8 & 34.40  & 37.40  & 13.98  & 8.02\% \\
rn.10v/100t.9 & 39.82  & 41.80  & 27.53  & 4.74\% \\
\midrule
Mean  & 40.56  & 42.94  & 19.06  & 5.62\% \\
\end{tabular}%

\begin{tablenotes}
\scriptsize
\item \textit{Note:} Aggregate gap to $Z_{MIP}$ is 5.55\%.
\end{tablenotes}
\end{threeparttable}
\end{table}
\begin{table}[!htbp]
\centering
\caption{Instance rn.11v/110t.}
\label{tab:11v/110t}
\scriptsize
\renewcommand{\arraystretch}{0.8}
\begin{threeparttable}
\begin{tabular}{ccccc}
\toprule 
\makecell[c]{Instance\\(idx)} & $Z_{\text{S-ALNS}}$ & $Z_{MIP}$ & Time (s) & \makecell[c]{Gap to\\$Z_{MIP}$} \\
\midrule
rn.11v/110t.0 & 53.20  & 53.04  & 54.57  & -0.30\% \\
rn.11v/110t.1 & 39.12  & 40.64  & 12.38  & 3.74\% \\
rn.11v/110t.2 & 47.20  & 49.76  & 56.27  & 5.14\% \\
rn.11v/110t.3 & 49.24  & 50.68  & 93.55  & 2.84\% \\
rn.11v/110t.4 & 46.60  & 44.74  & 25.77  & -4.16\% \\
rn.11v/110t.5 & 42.90  & 40.66  & 20.80  & -5.51\% \\
rn.11v/110t.6 & 40.04  & 45.16  & 19.74  & 11.34\% \\
rn.11v/110t.7 & 48.78  & 47.86  & 19.60  & -1.92\% \\
rn.11v/110t.8 & 42.26  & 44.12  & 17.96  & 4.22\% \\
rn.11v/110t.9 & 40.86  & 41.34  & 23.90  & 1.16\% \\
\midrule
Mean  & 45.02  & 45.80  & 34.45  & 1.66\% \\
\end{tabular}%

\begin{tablenotes}
\scriptsize
\item \textit{Note:} Aggregate gap to $Z_{MIP}$ is 1.70\%.
\end{tablenotes}
\end{threeparttable}
\end{table}
\begin{table}[!htbp]
\centering
\caption{Instance rn.13v/130t.}
\label{tab:13v/130t}
\scriptsize
\renewcommand{\arraystretch}{0.8}
\begin{threeparttable}
\begin{tabular}{ccccc}
\toprule 
\makecell[c]{Instance\\(idx)} & $Z_{\text{S-ALNS}}$ & $Z_{MIP}$ & Time (s) & \makecell[c]{Gap to\\$Z_{MIP}$} \\
\midrule
rn.13v/130t.0 & 59.78  & 58.76  & 46.67  & -1.74\% \\
rn.13v/130t.1 & 50.70  & 49.42  & 56.07  & -2.59\% \\
rn.13v/130t.2 & 49.08  & 50.50  & 28.99  & 2.81\% \\
rn.13v/130t.3 & 48.88  & 48.14  & 25.99  & -1.54\% \\
rn.13v/130t.4 & 50.94  & 45.28  & 27.87  & -12.50\% \\
rn.13v/130t.5 & 52.72  & 48.70  & 44.93  & -8.25\% \\
rn.13v/130t.6 & 56.00  & 47.44  & 40.78  & -18.04\% \\
rn.13v/130t.7 & 45.56  & 42.90  & 22.54  & -6.20\% \\
rn.13v/130t.8 & 47.62  & 45.14  & 28.82  & -5.49\% \\
rn.13v/130t.9 & 52.30  & 51.88  & 41.66  & -0.81\% \\
\midrule
Mean  & 51.36  & 48.82  & 36.43  & -5.44\% \\
\end{tabular}%

\begin{tablenotes}
\scriptsize
\item \textit{Note:} Aggregate gap to $Z_{MIP}$ is -5.21\%.
\end{tablenotes}
\end{threeparttable}
\end{table}
\begin{table}[!htbp]
\centering
\caption{Instance rn.15v/150t.}
\label{tab:15v/150t}
\scriptsize
\renewcommand{\arraystretch}{0.8}
\begin{threeparttable}
\begin{tabular}{ccccc}
\toprule 
\makecell[c]{Instance\\(idx)} & $Z_{\text{S-ALNS}}$ & $Z_{MIP}$ & Time (s) & \makecell[c]{Gap to\\$Z_{MIP}$} \\
\midrule
rn.15v/150t.0 & 71.86  & 65.72  & 77.48  & -9.34\% \\
rn.15v/150t.1 & 60.74  & 55.08  & 38.30  & -10.28\% \\
rn.15v/150t.2 & 64.82  & 47.10  & 70.17  & -37.62\% \\
rn.15v/150t.3 & 67.34  & 56.06  & 48.91  & -20.12\% \\
rn.15v/150t.4 & 62.98  & 47.98  & 61.80  & -31.26\% \\
rn.15v/150t.5 & 61.96  & 47.08  & 43.01  & -31.61\% \\
rn.15v/150t.6 & 52.80  & 48.96  & 50.49  & -7.84\% \\
rn.15v/150t.7 & 63.84  & 54.58  & 66.67  & -16.97\% \\
rn.15v/150t.8 & 70.12  & 44.10  & 108.73  & -59.00\% \\
rn.15v/150t.9 & 57.34  & 45.90  & 34.84  & -24.92\% \\
\midrule
Mean  & 63.38  & 51.26  & 60.04  & -24.90\% \\
\end{tabular}%

\begin{tablenotes}
\scriptsize
\item \textit{Note:} Aggregate gap to $Z_{MIP}$ is -23.65\%.
\end{tablenotes}
\end{threeparttable}
\end{table}

\clearpage

\section{Detail result on additional road network instances}\label{append3}

\begin{table}[!htbp]
\centering
\caption{Instance rn.5v/75t.}
\label{tab:5v/75t}
\scriptsize
\renewcommand{\arraystretch}{0.8}
\begin{threeparttable}
\begin{tabular}{ccccc}
\toprule 
\makecell[c]{Instance\\(idx)} & $Z_{\text{S-ALNS}}$ & $Z_{MIP}$ & Time (s) & \makecell[c]{Gap to\\$Z_{MIP}$} \\
\midrule
rn.5v/75t.0 & 26.76  & 28.52  & 6.85  & 6.17\% \\
rn.5v/75t.1 & 22.60  & 25.90  & 3.56  & 12.74\% \\
rn.5v/75t.2 & 20.18  & 22.28  & 2.78  & 9.43\% \\
rn.5v/75t.3 & 26.96  & 28.54  & 6.19  & 5.54\% \\
rn.5v/75t.4 & 21.74  & 23.32  & 4.35  & 6.78\% \\
rn.5v/75t.5 & 25.24  & 26.48  & 3.65  & 4.68\% \\
rn.5v/75t.6 & 21.36  & 23.70  & 3.54  & 9.87\% \\
rn.5v/75t.7 & 23.54  & 25.22  & 2.76  & 6.66\% \\
rn.5v/75t.8 & 24.36  & 26.42  & 7.44  & 7.80\% \\
rn.5v/75t.9 & 25.16  & 26.76  & 7.90  & 5.98\% \\
\midrule
Mean  & 23.79  & 25.71  & 4.90  & 7.56\% \\
\end{tabular}%

\begin{tablenotes}
\scriptsize
\item \textit{Note:} Aggregate gap to $Z_{MIP}$ is 7.48\%.
\end{tablenotes}
\end{threeparttable}
\end{table}

\begin{table}[!htbp]
\centering
\caption{Instance rn.5v/100t.}
\label{tab:5v/100t}
\scriptsize
\renewcommand{\arraystretch}{0.8}
\begin{threeparttable}
\begin{tabular}{ccccc}
\toprule 
\makecell[c]{Instance\\(idx)} & $Z_{\text{S-ALNS}}$ & $Z_{MIP}$ & Time (s) & \makecell[c]{Gap to\\$Z_{MIP}$} \\
\midrule
rn.5v/100t.0 & 32.14  & 35.40  & 3.01  & 9.21\% \\
rn.5v/100t.1 & 24.94  & 29.90  & 3.01  & 16.59\% \\
rn.5v/100t.2 & 19.68  & 21.40  & 5.88  & 8.04\% \\
rn.5v/100t.3 & 30.94  & 34.44  & 1.54  & 10.16\% \\
rn.5v/100t.4 & 31.46  & 33.56  & 3.06  & 6.26\% \\
rn.5v/100t.5 & 29.82  & 28.12  & 9.26  & -6.05\% \\
rn.5v/100t.6 & 30.66  & 33.00  & 4.52  & 7.09\% \\
rn.5v/100t.7 & 31.68  & 35.22  & 1.81  & 10.05\% \\
rn.5v/100t.8 & 29.96  & 20.38  & 2.42  & -47.01\% \\
rn.5v/100t.9 & 34.32  & 36.28  & 5.94  & 5.40\% \\
\midrule
Mean  & 29.56  & 30.77  & 4.05  & 1.97\% \\
\end{tabular}%

\begin{tablenotes}
\scriptsize
\item \textit{Note:} Aggregate gap to $Z_{MIP}$ is 3.93\%.
\end{tablenotes}
\end{threeparttable}
\end{table}

\begin{table}[!htbp]
\centering
\caption{Instance rn.7v/105t.}
\label{tab:7v/105t}
\scriptsize
\renewcommand{\arraystretch}{0.8}
\begin{threeparttable}
\begin{tabular}{ccccc}
\toprule 
\makecell[c]{Instance\\(idx)} & $Z_{\text{S-ALNS}}$ & $Z_{MIP}$ & Time (s) & \makecell[c]{Gap to\\$Z_{MIP}$} \\
\midrule
rn.7v/105t.0 & 29.44  & 33.34  & 9.55   & 11.70\% \\
rn.7v/105t.1 & 37.04  & 40.60  & 9.68   & 8.77\% \\
rn.7v/105t.2 & 33.20  & 35.84  & 14.26  & 7.37\% \\
rn.7v/105t.3 & 33.52  & 33.70  & 19.93  & 0.53\% \\
rn.7v/105t.4 & 37.24  & 37.34  & 30.94  & 0.27\% \\
rn.7v/105t.5 & 35.86  & 37.82  & 9.99   & 5.18\% \\
rn.7v/105t.6 & 30.90  & 32.58  & 8.52   & 5.16\% \\
rn.7v/105t.7 & 35.00  & 38.02  & 13.52  & 7.94\% \\
rn.7v/105t.8 & 33.74  & 35.02  & 4.99   & 3.66\% \\
rn.7v/105t.9 & 36.28  & 36.38  & 7.46   & 0.27\% \\
\midrule
Mean  & 34.22  & 36.06  & 12.88  & 5.08\% \\
\end{tabular}%

\begin{tablenotes}
\scriptsize
\item \textit{Note:} Aggregate gap to $Z_{MIP}$ is 5.11\%.
\end{tablenotes}
\end{threeparttable}
\end{table}

\begin{table}[!htbp]
\centering
\caption{Instance rn.7v/140t.}
\label{tab:7v/140t}
\scriptsize
\renewcommand{\arraystretch}{0.8}
\begin{threeparttable}
\begin{tabular}{ccccc}
\toprule 
\makecell[c]{Instance\\(idx)} & $Z_{\text{S-ALNS}}$ & $Z_{MIP}$ & Time (s) & \makecell[c]{Gap to\\$Z_{MIP}$} \\
\midrule
rn.7v/140t.0 & 32.20  & 33.02  & 8.12  & 2.48\% \\
rn.7v/140t.1 & 40.52  & 39.30  & 17.47  & -3.10\% \\
rn.7v/140t.2 & 35.08  & 37.68  & 10.11  & 6.90\% \\
rn.7v/140t.3 & 46.84  & 47.10  & 21.45  & 0.55\% \\
rn.7v/140t.4 & 36.56  & 37.86  & 9.19  & 3.43\% \\
rn.7v/140t.5 & 45.94  & 46.58  & 57.57  & 1.37\% \\
rn.7v/140t.6 & 36.06  & 34.90  & 9.51  & -3.32\% \\
rn.7v/140t.7 & 41.18  & 41.92  & 28.41  & 1.77\% \\
rn.7v/140t.8 & 47.14  & 48.76  & 16.62  & 3.32\% \\
rn.7v/140t.9 & 40.98  & 39.78  & 11.76  & -3.02\% \\
\midrule
Mean  & 40.25  & 40.69  & 19.02  & 1.04\% \\
\end{tabular}%

\begin{tablenotes}
\scriptsize
\item \textit{Note:} Aggregate gap to $Z_{MIP}$ is 1.08\%.
\end{tablenotes}
\end{threeparttable}
\end{table}

\begin{table}[!htbp]
\centering
\caption{Instance rn.9v/135t.}
\label{tab:9v/135t}
\scriptsize
\renewcommand{\arraystretch}{0.8}
\begin{threeparttable}
\begin{tabular}{ccccc}
\toprule 
\makecell[c]{Instance\\(idx)} & $Z_{\text{S-ALNS}}$ & $Z_{MIP}$ & Time (s) & \makecell[c]{Gap to\\$Z_{MIP}$} \\
\midrule
rn.9v/135t.0 & 35.38  & 36.28  & 13.61  & 2.48\% \\
rn.9v/135t.1 & 43.96  & 44.88  & 17.82  & 2.05\% \\
rn.9v/135t.2 & 35.72  & 40.86  & 7.89  & 12.58\% \\
rn.9v/135t.3 & 38.20  & 38.62  & 17.44  & 1.09\% \\
rn.9v/135t.4 & 48.46  & 45.50  & 24.17  & -6.51\% \\
rn.9v/135t.5 & 40.26  & 40.84  & 22.10  & 1.42\% \\
rn.9v/135t.6 & 38.24  & 43.10  & 11.99  & 11.28\% \\
rn.9v/135t.7 & 49.60  & 48.44  & 30.83  & -2.39\% \\
rn.9v/135t.8 & 42.66  & 42.12  & 15.47  & -1.28\% \\
rn.9v/135t.9 & 51.42  & 50.90  & 62.70  & -1.02\% \\
\midrule
Mean  & 42.39  & 43.15  & 22.40  & 1.97\% \\
\end{tabular}%

\begin{tablenotes}
\scriptsize
\item \textit{Note:} Aggregate gap to $Z_{MIP}$ is 1.77\%.
\end{tablenotes}
\end{threeparttable}
\end{table}

\begin{table}[!htbp]
\centering
\caption{Instance rn.9v/180t.}
\label{tab:9v/180t}
\scriptsize
\renewcommand{\arraystretch}{0.8}
\begin{threeparttable}
\begin{tabular}{ccccc}
\toprule 
\makecell[c]{Instance\\(idx)} & $Z_{\text{S-ALNS}}$ & $Z_{MIP}$ & Time (s) & \makecell[c]{Gap to\\$Z_{MIP}$} \\
\midrule
rn.9v/180t.0 & 60.78  & 46.68  & 20.01  & -30.21\% \\
rn.9v/180t.1 & 38.10  & 31.72  & 20.18  & -20.11\% \\
rn.9v/180t.2 & 52.36  & 44.34  & 24.10  & -18.09\% \\
rn.9v/180t.3 & 42.46  & 26.54  & 14.08  & -59.98\% \\
rn.9v/180t.4 & 54.18  & 44.32  & 24.43  & -22.25\% \\
rn.9v/180t.5 & 50.80  & 45.34  & 36.92  & -12.04\% \\
rn.9v/180t.6 & 42.76  & 37.06  & 16.35  & -15.38\% \\
rn.9v/180t.7 & 48.80  & 44.92  & 13.88  & -8.64\% \\
rn.9v/180t.8 & 41.48  & 19.06  & 15.08  & -117.63\% \\
rn.9v/180t.9 & 40.06  & 35.72  & 11.68  & -12.15\% \\
\midrule
Mean  & 47.18  & 37.57  & 19.67  & -31.65\% \\
\end{tabular}%

\begin{tablenotes}
\scriptsize
\item \textit{Note:} Aggregate gap to $Z_{MIP}$ is -25.57\%.
\end{tablenotes}
\end{threeparttable}
\end{table}

\begin{table}[!htbp]
\centering
\caption{Instance rn.10v/50t.}
\label{tab:10v/50t}
\scriptsize
\renewcommand{\arraystretch}{0.8}
\begin{threeparttable}
\begin{tabular}{ccccc}
\toprule 
\makecell[c]{Instance\\(idx)} & $Z_{\text{S-ALNS}}$ & $Z_{MIP}$ & Time (s) & \makecell[c]{Gap to\\$Z_{MIP}$} \\
\midrule
rn.10v/50t.0 & 27.28  & 27.82  & 3.99  & 1.94\% \\
rn.10v/50t.1 & 29.02  & 29.22  & 12.80  & 0.68\% \\
rn.10v/50t.2 & 27.92  & 28.06  & 11.51  & 0.50\% \\
rn.10v/50t.3 & 27.96  & 28.42  & 5.65  & 1.62\% \\
rn.10v/50t.4 & 28.76  & 29.00  & 20.98  & 0.83\% \\
rn.10v/50t.5 & 28.36  & 29.58  & 13.97  & 4.12\% \\
rn.10v/50t.6 & 30.56  & 30.48  & 9.05  & -0.26\% \\
rn.10v/50t.7 & 27.34  & 28.78  & 7.32  & 5.00\% \\
rn.10v/50t.8 & 20.92  & 21.88  & 4.12  & 4.39\% \\
rn.10v/50t.9 & 30.28  & 30.08  & 13.08  & -0.66\% \\
\midrule
Mean  & 27.84  & 28.33  & 10.25  & 1.82\% \\
\end{tabular}%

\begin{tablenotes}
\scriptsize
\item \textit{Note:} Aggregate gap to $Z_{MIP}$ is 1.74\%.
\end{tablenotes}
\end{threeparttable}
\end{table}

\begin{table}[!htbp]
\centering
\caption{Instance rn.10v/70t.}
\label{tab:10v/70t}
\scriptsize
\renewcommand{\arraystretch}{0.8}
\begin{threeparttable}
\begin{tabular}{ccccc}
\toprule 
\makecell[c]{Instance\\(idx)} & $Z_{\text{S-ALNS}}$ & $Z_{MIP}$ & Time (s) & \makecell[c]{Gap to\\$Z_{MIP}$} \\
\midrule
rn.10v/70t.0 & 33.52  & 34.10  & 10.29  & 1.70\% \\
rn.10v/70t.1 & 26.48  & 27.16  & 8.94  & 2.50\% \\
rn.10v/70t.2 & 35.06  & 36.34  & 10.90  & 3.52\% \\
rn.10v/70t.3 & 37.82  & 38.80  & 21.27  & 2.53\% \\
rn.10v/70t.4 & 35.80  & 35.06  & 21.81  & -2.11\% \\
rn.10v/70t.5 & 29.84  & 31.80  & 6.42  & 6.16\% \\
rn.10v/70t.6 & 28.46  & 29.52  & 6.61  & 3.59\% \\
rn.10v/70t.7 & 30.02  & 31.76  & 8.24  & 5.48\% \\
rn.10v/70t.8 & 37.46  & 37.82  & 19.04  & 0.95\% \\
rn.10v/70t.9 & 33.80  & 36.06  & 11.18  & 6.27\% \\
\midrule
Mean  & 32.83  & 33.84  & 12.47  & 3.06\% \\
\end{tabular}%

\begin{tablenotes}
\scriptsize
\item \textit{Note:} Aggregate gap to $Z_{MIP}$ is 3.00\%.
\end{tablenotes}
\end{threeparttable}
\end{table}

\begin{table}[!htbp]
\centering
\caption{Instance rn.10v/90t.}
\label{tab:10v/90t}
\scriptsize
\renewcommand{\arraystretch}{0.8}
\begin{threeparttable}
\begin{tabular}{ccccc}
\toprule 
\makecell[c]{Instance\\(idx)} & $Z_{\text{S-ALNS}}$ & $Z_{MIP}$ & Time (s) & \makecell[c]{Gap to\\$Z_{MIP}$} \\
\midrule
rn.10v/90t.0 & 39.90  & 42.34  & 11.84  & 5.76\% \\
rn.10v/90t.1 & 36.96  & 40.70  & 17.56  & 9.19\% \\
rn.10v/90t.2 & 36.46  & 37.74  & 9.24  & 3.39\% \\
rn.10v/90t.3 & 37.64  & 38.34  & 9.47  & 1.83\% \\
rn.10v/90t.4 & 41.30  & 43.02  & 8.49  & 4.00\% \\
rn.10v/90t.5 & 37.64  & 37.68  & 12.74  & 0.11\% \\
rn.10v/90t.6 & 43.36  & 44.84  & 37.26  & 3.30\% \\
rn.10v/90t.7 & 41.32  & 42.02  & 42.55  & 1.67\% \\
rn.10v/90t.8 & 40.88  & 44.60  & 13.18  & 8.34\% \\
rn.10v/90t.9 & 34.82  & 36.88  & 14.28  & 5.59\% \\
\midrule
Mean  & 39.03  & 40.82  & 17.66  & 4.32\% \\
\end{tabular}%

\begin{tablenotes}
\scriptsize
\item \textit{Note:} Aggregate gap to $Z_{MIP}$ is 4.38\%.
\end{tablenotes}
\end{threeparttable}
\end{table}

\begin{table}[!htbp]
\centering
\caption{Instance rn.10v/110t.}
\label{tab:10v/110t}
\scriptsize
\renewcommand{\arraystretch}{0.8}
\begin{threeparttable}
\begin{tabular}{ccccc}
\toprule 
\makecell[c]{Instance\\(idx)} & $Z_{\text{S-ALNS}}$ & $Z_{MIP}$ & Time (s) & \makecell[c]{Gap to\\$Z_{MIP}$} \\
\midrule
rn.10v/110t.0 & 41.30  & 43.48  & 45.81  & 5.01\% \\
rn.10v/110t.1 & 40.20  & 42.94  & 12.95  & 6.38\% \\
rn.10v/110t.2 & 41.90  & 43.04  & 20.51  & 2.65\% \\
rn.10v/110t.3 & 46.64  & 48.36  & 13.57  & 3.56\% \\
rn.10v/110t.4 & 47.66  & 47.60  & 58.06  & -0.13\% \\
rn.10v/110t.5 & 38.44  & 41.38  & 28.09  & 7.10\% \\
rn.10v/110t.6 & 49.16  & 48.90  & 22.43  & -0.53\% \\
rn.10v/110t.7 & 47.88  & 47.46  & 18.65  & -0.88\% \\
rn.10v/110t.8 & 42.18  & 43.62  & 27.61  & 3.30\% \\
rn.10v/110t.9 & 43.52  & 45.30  & 9.42  & 3.93\% \\
\midrule
Mean          & 43.89  & 45.21  & 25.71  & 3.04\% \\
\end{tabular}%

\begin{tablenotes}
\scriptsize
\item \textit{Note:} Aggregate gap to $Z_{MIP}$ is 2.92\%.
\end{tablenotes}
\end{threeparttable}
\end{table}

\begin{table}[!htbp]
\centering
\caption{Instance rn.10v/130t.}
\label{tab:10v/130t}
\scriptsize
\renewcommand{\arraystretch}{0.8}
\begin{threeparttable}
\begin{tabular}{ccccc}
\toprule 
\makecell[c]{Instance\\(idx)} & $Z_{\text{S-ALNS}}$ & $Z_{MIP}$ & Time (s) & \makecell[c]{Gap to\\$Z_{MIP}$} \\
\midrule
rn.10v/130t.0 & 43.58  & 45.86  & 12.89  & 4.97\% \\
rn.10v/130t.1 & 41.96  & 41.34  & 13.58  & -1.50\% \\
rn.10v/130t.2 & 58.88  & 59.56  & 64.51  & 1.14\% \\
rn.10v/130t.3 & 52.52  & 49.78  & 41.59  & -5.50\% \\
rn.10v/130t.4 & 52.52  & 51.16  & 27.95  & -2.66\% \\
rn.10v/130t.5 & 50.76  & 48.56  & 21.21  & -4.53\% \\
rn.10v/130t.6 & 45.48  & 48.94  & 28.35  & 7.07\% \\
rn.10v/130t.7 & 45.14  & 44.20  & 21.64  & -2.13\% \\
rn.10v/130t.8 & 48.22  & 46.54  & 30.98  & -3.61\% \\
rn.10v/130t.9 & 38.64  & 35.20  & 9.88  & -9.77\% \\
\midrule
Mean          & 47.77  & 47.11  & 27.26  & -1.65\% \\
\end{tabular}%

\begin{tablenotes}
\scriptsize
\item \textit{Note:} Aggregate gap to $Z_{MIP}$ is -1.39\%.
\end{tablenotes}
\end{threeparttable}
\end{table}

\begin{table}[!htbp]
\centering
\caption{Instance rn.10v/150t.}
\label{tab:10v/150t}
\scriptsize
\renewcommand{\arraystretch}{0.8}
\begin{threeparttable}
\begin{tabular}{ccccc}
\toprule 
\makecell[c]{Instance\\(idx)} & $Z_{\text{S-ALNS}}$ & $Z_{MIP}$ & Time (s) & \makecell[c]{Gap to\\$Z_{MIP}$} \\
\midrule
rn.10v/150t.0 & 53.14  & 47.84  & 24.48  & -11.08\% \\
rn.10v/150t.1 & 48.44  & 47.46  & 30.34  & -2.06\% \\
rn.10v/150t.2 & 52.92  & 52.70  & 25.84  & -0.42\% \\
rn.10v/150t.3 & 45.92  & 44.18  & 26.93  & -3.94\% \\
rn.10v/150t.4 & 46.80  & 44.56  & 34.64  & -5.03\% \\
rn.10v/150t.5 & 49.90  & 46.30  & 23.87  & -7.78\% \\
rn.10v/150t.6 & 53.16  & 49.68  & 18.59  & -7.00\% \\
rn.10v/150t.7 & 39.22  & 40.78  & 11.87  & 3.83\% \\
rn.10v/150t.8 & 48.98  & 46.94  & 22.25  & -4.35\% \\
rn.10v/150t.9 & 53.54  & 52.82  & 37.29  & -1.36\% \\
\midrule
Mean          & 49.20  & 47.33  & 25.61  & -3.92\% \\
\end{tabular}%

\begin{tablenotes}
\scriptsize
\item \textit{Note:} Aggregate gap to $Z_{MIP}$ is -3.96\%.
\end{tablenotes}
\end{threeparttable}
\end{table}

\begin{table}[!htbp]
\centering
\caption{Instance rn.11v/165t.}
\label{tab:11v/165t}
\scriptsize
\renewcommand{\arraystretch}{0.8}
\begin{threeparttable}
\begin{tabular}{ccccc}
\toprule 
\makecell[c]{Instance\\(idx)} & $Z_{\text{S-ALNS}}$ & $Z_{MIP}$ & Time (s) & \makecell[c]{Gap to\\$Z_{MIP}$} \\
\midrule
rn.11v/165t.0 & 48.46  & 33.80  & 17.41  & -43.37\% \\
rn.11v/165t.1 & 58.78  & 38.90  & 27.72  & -51.11\% \\
rn.11v/165t.2 & 57.52  & 48.50  & 27.68  & -18.60\% \\
rn.11v/165t.3 & 64.34  & 40.02  & 99.31  & -60.77\% \\
rn.11v/165t.4 & 54.94  & 44.12  & 30.57  & -24.52\% \\
rn.11v/165t.5 & 55.10  & 35.52  & 29.37  & -55.12\% \\
rn.11v/165t.6 & 55.40  & 39.88  & 42.76  & -38.92\% \\
rn.11v/165t.7 & 58.00  & 45.44  & 29.76  & -27.64\% \\
rn.11v/165t.8 & 63.40  & 42.44  & 58.29  & -49.39\% \\
rn.11v/165t.9 & 50.66  & 36.60  & 25.81  & -38.42\% \\
\midrule
Mean  & 56.66  & 40.52  & 38.87  & -40.79\% \\
\end{tabular}%

\begin{tablenotes}
\scriptsize
\item \textit{Note:} Aggregate gap to $Z_{MIP}$ is -39.83\%.
\end{tablenotes}
\end{threeparttable}
\end{table}

\begin{table}[!htbp]
\centering
\caption{Instance rn.11v/220t.}
\label{tab:11v/220t}
\scriptsize
\renewcommand{\arraystretch}{0.8}
\begin{threeparttable}
\begin{tabular}{ccccc}
\toprule 
\makecell[c]{Instance\\(idx)} & $Z_{\text{S-ALNS}}$ & $Z_{MIP}$ & Time (s) & \makecell[c]{Gap to\\$Z_{MIP}$} \\
\midrule
rn.11v/220t.0 & 61.88  & 42.08  & 33.60  & -47.05\% \\
rn.11v/220t.1 & 62.34  & 36.84  & 48.42  & -69.22\% \\
rn.11v/220t.2 & 56.38  & 30.28  & 42.91  & -86.20\% \\
rn.11v/220t.3 & 65.58  & 42.32  & 72.04  & -54.96\% \\
rn.11v/220t.4 & 78.34  & 52.34  & 122.78  & -49.68\% \\
rn.11v/220t.5 & 58.20  & 37.20  & 27.91  & -56.45\% \\
rn.11v/220t.6 & 53.08  & 29.92  & 44.91  & -77.41\% \\
rn.11v/220t.7 & 62.62  & 40.84  & 87.72  & -53.33\% \\
rn.11v/220t.8 & 50.94  & 28.62  & 17.75  & -77.99\% \\
rn.11v/220t.9 & 55.38  & 40.06  & 40.17  & -38.24\% \\
\midrule
Mean  & 60.47  & 38.05  & 53.82  & -61.05\% \\
\end{tabular}%

\begin{tablenotes}
\scriptsize
\item \textit{Note:} Aggregate gap to $Z_{MIP}$ is -58.93\%.
\end{tablenotes}
\end{threeparttable}
\end{table}

\begin{table}[!htbp]
\centering
\caption{Instance rn.13v/195t.}
\label{tab:13v/195t}
\scriptsize
\renewcommand{\arraystretch}{0.8}
\begin{threeparttable}
\begin{tabular}{ccccc}
\toprule 
\makecell[c]{Instance\\(idx)} & $Z_{\text{S-ALNS}}$ & $Z_{MIP}$ & Time (s) & \makecell[c]{Gap to\\$Z_{MIP}$} \\
\midrule
rn.13v/195t.0 & 62.68  & 39.58  & 46.76  & -58.36\% \\
rn.13v/195t.1 & 73.56  & 52.94  & 143.75  & -38.95\% \\
rn.13v/195t.2 & 63.46  & 41.64  & 169.78  & -52.40\% \\
rn.13v/195t.3 & 66.46  & 35.00  & 84.11  & -89.89\% \\
rn.13v/195t.4 & 64.60  & 44.12  & 30.10  & -46.42\% \\
rn.13v/195t.5 & 57.94  & 42.14  & 39.84  & -37.49\% \\
rn.13v/195t.6 & 66.72  & 47.60  & 61.99  & -40.17\% \\
rn.13v/195t.7 & 69.48  & 32.66  & 48.35  & -112.74\% \\
rn.13v/195t.8 & 59.34  & 32.78  & 39.53  & -81.03\% \\
rn.13v/195t.9 & 73.72  & 39.22  & 90.59  & -87.97\% \\
\midrule
Mean  & 65.80  & 40.77  & 75.48  & -64.54\% \\
\end{tabular}%

\begin{tablenotes}
\scriptsize
\item \textit{Note:} Aggregate gap to $Z_{MIP}$ is -61.39\%.
\end{tablenotes}
\end{threeparttable}
\end{table}

\begin{table}[!htbp]
\centering
\caption{Instance rn.13v/260t.}
\label{tab:13v/260t}
\scriptsize
\renewcommand{\arraystretch}{0.8}
\begin{threeparttable}
\begin{tabular}{ccccc}
\toprule 
\makecell[c]{Instance\\(idx)} & $Z_{\text{S-ALNS}}$ & $Z_{MIP}$ & Time (s) & \makecell[c]{Gap to\\$Z_{MIP}$} \\
\midrule
rn.13v/260t.0 & 63.62  & 29.82  & 41.62  & -113.35\% \\
rn.13v/260t.1 & 75.34  & 34.06  & 169.50  & -121.20\% \\
rn.13v/260t.2 & 71.60  & 30.54  & 55.53  & -134.45\% \\
rn.13v/260t.3 & 65.54  & 27.02  & 63.61  & -142.56\% \\
rn.13v/260t.4 & 64.66  & 26.84  & 52.83  & -140.91\% \\
rn.13v/260t.5 & 78.42  & 26.62  & 236.02  & -194.59\% \\
rn.13v/260t.6 & 70.54  & 34.50  & 62.19  & -104.46\% \\
rn.13v/260t.7 & 78.00  & 12.50  & 63.89  & -524.00\% \\
rn.13v/260t.8 & 68.16  & 32.12  & 39.36  & -112.20\% \\
rn.13v/260t.9 & 76.58  & 41.16  & 159.85  & -86.05\% \\
\midrule
Mean  & 71.25  & 29.52  & 94.44  & -167.38\% \\
\end{tabular}%

\begin{tablenotes}
\scriptsize
\item \textit{Note:} Aggregate gap to $Z_{MIP}$ is -141.36\%.
\end{tablenotes}
\end{threeparttable}
\end{table}

\begin{table}[!htbp]
\centering
\caption{Instance rn.15v/225t.}
\label{tab:15v/225t}
\scriptsize
\renewcommand{\arraystretch}{0.8}
\begin{threeparttable}
\begin{tabular}{ccccc}
\toprule 
\makecell[c]{Instance\\(idx)} & $Z_{\text{S-ALNS}}$ & $Z_{MIP}$ & Time (s) & \makecell[c]{Gap to\\$Z_{MIP}$} \\
\midrule
rn.15v/225t.0 & 72.18  & 37.82  & 82.29  & -90.85\% \\
rn.15v/225t.1 & 69.40  & 33.84  & 60.26  & -105.08\% \\
rn.15v/225t.2 & 73.84  & 27.78  & 77.29  & -165.80\% \\
rn.15v/225t.3 & 75.82  & 41.10  & 97.32  & -84.48\% \\
rn.15v/225t.4 & 58.68  & 28.02  & 32.16  & -109.42\% \\
rn.15v/225t.5 & 68.10  & 31.88  & 52.43  & -113.61\% \\
rn.15v/225t.6 & 68.84  & 14.68  & 51.97  & -368.94\% \\
rn.15v/225t.7 & 57.36  & 24.62  & 24.79  & -132.98\% \\
rn.15v/225t.8 & 70.38  & 35.82  & 66.97  & -96.48\% \\
rn.15v/225t.9 & 78.18  & 40.84  & 74.14  & -91.43\% \\
\midrule
Mean  & 69.28  & 31.64  & 61.96  & -135.91\% \\
\end{tabular}%

\begin{tablenotes}
\scriptsize
\item \textit{Note:} Aggregate gap to $Z_{MIP}$ is -118.96\%.
\end{tablenotes}
\end{threeparttable}
\end{table}

\begin{table}[!htbp]
\centering
\caption{Instance rn.15v/300t.}
\label{tab:15v/300t}
\scriptsize
\renewcommand{\arraystretch}{0.8}
\begin{threeparttable}
\begin{tabular}{ccccc}
\toprule 
\makecell[c]{Instance\\(idx)} & $Z_{\text{S-ALNS}}$ & $Z_{MIP}$ & Time (s) & \makecell[c]{Gap to\\$Z_{MIP}$} \\
\midrule
rn.15v/300t.0 & 84.98  & 0.00   & 96.95    & - \\
rn.15v/300t.1 & 79.50  & 26.68  & 133.52  & -197.98\% \\
rn.15v/300t.2 & 84.66  & 34.98  & 128.23  & -142.02\% \\
rn.15v/300t.3 & 85.02  & 14.08  & 120.21  & -503.84\% \\
rn.15v/300t.4 & 80.54  & 23.62  & 63.32   & -240.98\% \\
rn.15v/300t.5 & 85.88  & 21.58  & 144.16  & -297.96\% \\
rn.15v/300t.6 & 83.96  & 34.38  & 88.90   & -144.21\% \\
rn.15v/300t.7 & 90.00  & 27.20  & 151.58  & -230.88\% \\
rn.15v/300t.8 & 81.98  & 22.04  & 152.90  & -271.96\% \\
rn.15v/300t.9 & 70.20  & 18.32  & 40.62   & -283.19\% \\
\midrule
Mean  & 82.67  & 22.29  & 112.04  & -257.00\% \\
\end{tabular}%

\begin{tablenotes}
\scriptsize
\item \textit{Note:} Aggregate gap to $Z_{MIP}$ is -270.93\%.
\end{tablenotes}
\end{threeparttable}
\end{table}

\end{document}